\documentclass[11pt,reqno]{amsart}

\usepackage[a4paper,left=35mm,right=35mm,top=30mm,bottom=30mm,marginpar=25mm]{geometry}

\usepackage{amsmath}
\usepackage{amssymb}
\usepackage{amsthm}
\usepackage[utf8]{inputenc}
\usepackage{graphicx}
\usepackage{placeins}
\usepackage{hyperref}
\usepackage{dsfont}
\usepackage[nocompress, space]{cite}
\usepackage{xcolor}

\usepackage[displaymath,mathlines]{lineno}
\allowdisplaybreaks

\newcommand{\tr}{\operatorname{tr}}

\newcommand{\supp}{\operatorname{supp}}

\renewcommand{\div}{\operatorname{div}}

\newcommand{\diag}{\operatorname{diag}}

\newcommand{\Rr}{{\mathbb{R}}}

\newcommand{\Nn}{{\mathbb{N}}}

\newcommand{\cws}{{\overset{*}{\rightharpoonup}}}

\def\leq{\leqslant}
\def\geq{\geqslant}

\numberwithin{equation}{section}

\newtheoremstyle{thmlemcorr}{10pt}{10pt}{\itshape}{}{\bfseries}{.}{10pt}{{\thmname{#1}\thmnumber{
#2}\thmnote{ (#3)}}}
\newtheoremstyle{thmlemcorr*}{10pt}{10pt}{\itshape}{}{\bfseries}{.}\newline{{\thmname{#1}\thmnumber{
#2}\thmnote{ (#3)}}}
\newtheoremstyle{defi}{10pt}{10pt}{\itshape}{}{\bfseries}{.}{10pt}{{\thmname{#1}\thmnumber{
#2}\thmnote{ (#3)}}}
\newtheoremstyle{remexample}{10pt}{10pt}{}{}{\bfseries}{.}{10pt}{{\thmname{#1}\thmnumber{
#2}\thmnote{ (#3)}}}
\newtheoremstyle{ass}{10pt}{10pt}{}{}{\bfseries}{.}{10pt}{{\thmname{#1}\thmnumber{
A#2}\thmnote{ (#3)}}}

\theoremstyle{thmlemcorr}
\newtheorem{theorem}{Theorem}
\numberwithin{theorem}{section}
\newtheorem{lemma}[theorem]{Lemma}
\newtheorem{corollary}[theorem]{Corollary}
\newtheorem{proposition}[theorem]{Proposition}

\theoremstyle{thmlemcorr*}
\newtheorem{theorem*}{Theorem}
\newtheorem{lemma*}[theorem]{Lemma}
\newtheorem{corollary*}[theorem]{Corollary}
\newtheorem{proposition*}[theorem]{Proposition}
\newtheorem{problem*}[theorem]{Problem}
\newtheorem{conjecture*}[theorem]{Conjecture}

\theoremstyle{defi}

\newtheorem{hyp}{Assumption}

\theoremstyle{remexample}
\newtheorem{remark}[theorem]{Remark}
\newtheorem{example}[theorem]{Example}

\theoremstyle{ass}

\title[Collision intensity minimization in crowd flow]{Quantifying and minimizing collision intensity in stationary crowd flow}

\author{Yuri Ashrafyan}
\address{CEMSE Division, King Abdullah University of Science and Technology (KAUST), Thuwal 23955--6900, Saudi Arabia}
\email{yuri.ashrafyan@kaust.edu.sa}

\author{Diogo A. Gomes}
\address{CEMSE Division, King Abdullah University of Science and Technology (KAUST), Thuwal 23955--6900, Saudi Arabia}
\email{diogo.gomes@kaust.edu.sa}

\thanks{
The authors were partially supported by KAUST Baseline funds. 
}
\date{\today}

\begin{document}


\begin{abstract}
We develop a continuum model to quantify kinematic collision risk
in stationary crowd flow. From the density-weighted inward relative
flux, we derive a collision functional based on the
symmetric velocity gradient. For incompressible flows
with constant density, this functional is convex and positively
homogeneous of degree one in a combination of second derivatives
of the stream function. We minimize it subject to prescribed
boundary flux. On bounded simply connected domains with Lipschitz
boundary, under a finite-energy lifting assumption, we prove
existence of minimizers in a relaxed class that allows this
derivative combination to be a vector-valued Radon measure.
We show that disks are the only domains with a nontrivial
zero-trace kernel, corresponding to rigid rotation about the
center. Using the compact embedding for functions of bounded
deformation, we prove that relaxed energy sublevels are compact
without normalization on every other admissible domain, and after
removing one rotational component on disks. We identify the
relaxed functional with 
the strong \(W^{1,1}\) lower semicontinuous envelope of the functional on the prescribed-trace class with integrable controlled derivatives.
Then, we derive optimality conditions via
convex duality. We construct explicit radial flows on annular sectors and prove
that they minimize the energy over the full admissible classes,
including nonradial competitors. On the square, distinct
minimizers with the same prescribed boundary flux have different
local collision intensities. The square--disk comparison shows
that the saddle flow loses optimality when the domain changes,
even after normalization.
\end{abstract}

\maketitle

\section{Introduction}\label{sec:introduction}
Models of pedestrian dynamics have applications in urban planning, event management, and public safety. A key goal is to understand and optimize stationary crowd flows---persistent, large-scale movement patterns that are statistically constant over time. Such flows are described by a macroscopic density $m(x)$ and a velocity field $v(x)$, governed by the continuity equation
\[
    \div(mv)=0,
\]
which expresses conservation of agents.

Crowd flow models include the Hughes continuum framework
\cite{hughes2002continuum} and its analytical and numerical developments
such as in \cite{Gomes:2016aa} and references therein.
Mean-field-game formulations account for anticipation, crowd aversion,
and congestion \cite{Lachapelle10}, while optimal-control formulations
incorporate nonlinear mobility in evacuation problems
\cite{burgerMeanFieldGames2014}.
Non-separable mean-field games extend the Hughes framework to describe
movement through or around congested regions
\cite{ghattassiNonseparableMeanField2025}.
For dense crowds, density-constrained models restrict motion through
admissible-velocity projections \cite{mauryMacroscopicCrowdMotion2010},
while hybrid agent--continuum simulations use unilateral
incompressibility to model crowding and facilitate collision avoidance
\cite{narain2009aggregate}.
A related granular approach obtains congested crowd motion as the
$p\to\infty$ limit of Hele--Shaw-type problems with $p$-Laplacian
diffusion and drift \cite{IU24}.

Variational formulations for first-order stationary mean-field games
with congestion provide a basis for existence results and numerical
approximation \cite{EFGNV2017}.
For stationary crowd entry--exit dynamics on bounded domains,
mixed boundary conditions prescribe inflow and regulate exit through
an exit-cost complementarity condition,
while a variational approach yields existence of solutions
\cite{alharbi2026first}.
Prescribed boundary inflow has also been studied in time-dependent
first-order mean-field games through a relaxed variational formulation
and Fenchel--Rockafellar duality
\cite{gomesTimeDependentFirstOrder2024}.
In these formulations, costs depend on density and individual motion
rather than directly on the symmetric velocity gradient.
More generally, continuity and incompressibility constraints alone
do not penalize local shear or directional compression.

Anticipatory pedestrian models use projected time to collision
\cite{karamouzasUniversalPowerLaw2014}, while vision-based kinetic
and macroscopic models use bearing-angle changes and time to
interaction \cite{DARPT2013}.
Kinetic models also describe collision avoidance through sidestepping
\cite{festa2016kinetic}, while microscopic and macroscopic contact
models estimate the forces and pressure generated in dense crowds
\cite{jebrane2019estimating}.
In the deterministic kinematic interpretation adopted here,
agents in a smooth stationary velocity field follow the streamlines
$\dot{x}(t)=v(x(t))$, and collision risk refers to the instantaneous
approach of nearby agents, rather than intersections of point
trajectories.
Such approach depends on the local directional compression and shear
encoded in the velocity gradient and can occur even at moderate density.
At the same density, rigid rotation produces no local approach,
while shear and directional compression can do so.
This motivates using the symmetric velocity gradient to compare stationary
flow geometries under prescribed boundary flux.

We introduce a collision functional and develop a
variational framework for its minimization, addressing
two questions: how should collision intensity be quantified from the pair $(m,v)$, and, in the constant-density regime under prescribed boundary flux, which incompressible flow geometries minimize it?
Our derivation of the local collision intensity adapts the classical
Saffman--Turner inward-relative-flux construction for passive particles
in fluid flows \cite{saffmanCollisionDropsTurbulent1956} and its two-dimensional
eigenvalue evaluation in \cite[Appendix~A.1]{gustavssonCollisionsParticlesAdvected2008}.
In the aerosol and turbulence literature, this construction is used to
estimate collision rates of finite-radius particles advected in random
velocity fields. The Saffman--Turner approximation describes an initial
transient; long-time first-contact rates must also account for trajectory
re-entry and, in compressible flows, clustering
\cite{gustavssonCollisionsParticlesAdvected2008}.
Here, by contrast, we pass to the small-radius limit to define a pointwise
encounter-intensity density for a deterministic stationary flow on a
bounded domain, and study the variational boundary-value problem of
minimizing the resulting total collision functional under prescribed
boundary flux.

While the collision functional is derived for a general density $m$
and velocity field $v$, our variational existence and relaxation theory
focuses on constant-density incompressible flows ($m\equiv1$, $\div v=0$)
on bounded simply connected Lipschitz planar domains.
By introducing a stream function $\psi$ such that $v=D^\perp \psi$, we show that the dense collision functional becomes, up to the constant factor $1/\pi$, the convex functional
$$
\psi\mapsto\int_\Omega\sqrt{4\psi_{xy}^2+(\psi_{xx}-\psi_{yy})^2}\,dx\,dy,
$$
to be minimized over stream functions with the boundary trace prescribed by the flux.
Under a compatible boundary flux and a finite-energy lifting condition (Assumptions~\ref{ass:dense_boundary}
and~\ref{ass:dense_lifting}), we establish the existence of
relaxed minimizers and identify the relaxed functional with the lower semicontinuous
envelope of the unrelaxed functional in the strong
\(W^{1,1}\) topology. Using convex duality and integration by parts, we also obtain a dual representation
of the minimum, weak Euler--Lagrange conditions, and a smooth sufficient optimality criterion involving a natural boundary condition. We apply these conditions to prove that radial flow is a global minimizer in an annular sector, exhibit a family of minimizers on the square with different peak local collision intensities, and show how domain geometry rules out optimality of the saddle flow on a disk.

The proof of existence relies on the direct method in the calculus
of variations. The energy controls
\(TD^2\psi=(2\psi_{xy},\psi_{xx}-\psi_{yy})\), but not the
\(L^1\) norm of the full Hessian, and is therefore not coercive
in \(W^{2,1}(\Omega)\); see \cite{ornsteinNoninequalityDifferentialOperators1962}.
Since the controlled derivatives may concentrate along
bounded-energy sequences, we relax the problem by allowing
\(TD^2\psi\) to be a vector-valued Radon measure.

The prescribed scalar trace removes every zero-energy perturbation
except on disks, where the zero-trace kernel is one-dimensional
and corresponds to rigid rotation about the center. We prove
\(W^{1,1}(\Omega)\)-coercivity after subtracting the projection
onto this kernel. This preserves the prescribed trace and the
entire measure \(TD^2\psi\), although it may change the velocity.
Using the compact embedding of the space of functions of bounded deformation,
\(BD(\Omega)\hookrightarrow\hookrightarrow L^1(\Omega;\mathbb R^2)\),
recalled in \cite[Theorem~2.10]{barrosoRelaxationTheoremSpace2000}, we show that
relaxed energy sublevels are compact in strong \(W^{1,1}(\Omega)\)
without normalization on domains other than disks, and after
normalization on disks. Related relaxation problems for
linear-growth symmetric-gradient energies are studied in
\cite{barrosoRelaxationTheoremSpace2000}.
The dense collision energy is a constant multiple of the total variation of the symmetric velocity gradient; our analysis incorporates incompressibility, prescribed normal flux, and the resulting geometric kernel.

For exact relaxation, we use interior
mollification with a locally finite partition of unity. The smoothing
scales are chosen so that the difference from the original stream
function belongs to \(W^{1,1}_0(\Omega)\), while the collision energy
converges. Thus the prescribed trace is preserved exactly, and
incompressibility holds automatically throughout the approximation.

The paper is organized as follows.
Section~\ref{sec:notation} fixes notation.
Section~\ref{sec:collision_functional} derives the
collision-intensity functional and expresses it in terms of the
density and the eigenvalues of the symmetric velocity gradient.
Section~\ref{sec:dense_flow} applies the model to dense incompressible flows, introduces the boundary flux and lifting assumptions, and formulates the stream-function variational problem. 
Section~\ref{sec:relaxed_dense_problem} proves existence of relaxed minimizers, identifies the relaxed problem with the lower semicontinuous envelope of the unrelaxed one, and derives duality and optimality conditions.
Section~\ref{sec:Examples} studies radial flow in an annular region, characterizes zero-energy minimizers, exhibits a one-parameter family of minimizers on the square, and compares the saddle flow on a disk.

\section{Notation}\label{sec:notation}

Unless otherwise specified, $\Omega\subset\mathbb{R}^2$ denotes a
bounded simply connected domain with Lipschitz boundary
$\Gamma=\partial\Omega$, and $ds$ denotes arclength measure.
Along $\Gamma$, we denote the outer unit normal by $\nu=(\nu_1,\nu_2)$
and the positively oriented unit tangent vector by
\[
\tau:=(-\nu_2,\nu_1).
\]

For a scalar function $f$ smooth up to $\Gamma$, we write
$Df=(\partial_x f,\partial_y f)$, set
\[
D^\perp f:= (\partial_y f,-\partial_x f),
\]
and define the tangential derivative of $f$ along $\Gamma$ by
\[
\partial_\tau f := D f\cdot\tau = D^\perp f\cdot\nu.
\]
For boundary functions in \(W^{1,1}(\Gamma)\), tangential
derivatives are understood weakly with respect to arclength.
The positive and negative parts of a real-valued quantity $a$ are defined by
\[
a^+:=\max\{a,0\},\qquad a^-:=\max\{-a,0\}.
\]

Let $B_r(x)$ denote the open ball of radius $r>0$
centered at $x\in\mathbb R^2$.
Its boundary $\partial B_r(x)$ has outer unit normal $\nu(y)=(y-x)/r$, and we parametrize $\partial B_r(x)$ by $y=x+r\,e(\theta)$, where
\begin{equation}\label{etheta}
e(\theta):=(\cos\theta,\sin\theta), \qquad \theta\in[0,2\pi).
\end{equation}

We use standard notation for Lebesgue and Sobolev spaces, $L^p(\Omega)$ and $W^{k,p}(\Omega)$.
On Lipschitz boundaries, traces of $W^{1,1}(\Omega)$ functions (denoted simply by $\psi$ on $\Gamma$) are understood in the $L^1(\Gamma)$ sense, and $W^{1,1}_0(\Omega)$ denotes the subspace of functions in $W^{1,1}(\Omega)$ with zero trace.
We use $\mathcal M(\Omega;\mathbb R^k)$ for the space of finite
$\mathbb R^k$-valued Radon measures on $\Omega$. For $\mu\in\mathcal M(\Omega;\mathbb R^k)$,
$|\mu|$ denotes its total variation measure and
$|\mu|(\Omega)$ its total mass.
Weak-$*$ convergence in $\mathcal M(\Omega;\mathbb R^k)$, denoted by $\cws$, is understood in duality with
$C_0(\Omega;\mathbb R^k)$.
We identify an integrable field $f$ with the measure
$f\,\mathcal L^2$, where $\mathcal L^2$ denotes planar Lebesgue
measure.

For $w\in L^1(\Omega;\mathbb R^2)$, let
$Ew=(Dw+Dw^\top)/2$ in the sense of distributions, and set
\[
BD(\Omega):=
\left\{w\in L^1(\Omega;\mathbb R^2):
Ew\in\mathcal M(\Omega;\mathbb R_{\mathrm{sym}}^{2\times2})\right\}.
\]
Here $\mathbb R_{\mathrm{sym}}^{2\times2}$ carries the Frobenius norm,
and $\|w\|_{BD(\Omega)}:=\|w\|_{L^1(\Omega)}+|Ew|(\Omega)$.

Throughout, generic positive constants may change from line to line and depend
only on the fixed problem data.

\section{The collision functional}
\label{sec:collision_functional}

Throughout this section, we work under the assumption:

\begin{hyp}[Local regularity]\label{ass:local_regularity}
The velocity field and density satisfy
\[
    v\in C^{1,\alpha}_{\mathrm{loc}}(\Omega;\mathbb R^2),
    \qquad
    m\in C^{0,1}_{\mathrm{loc}}(\Omega),
    \qquad
    m\ge0,
\]
for some \(\alpha\in(0,1]\).
\end{hyp}

We quantify collision risk through the instantaneous approach of nearby agents, rather than intersections of point trajectories. To isolate this effect from uniform drift, we measure the influx of agents into a small ball $B_r(x)$, $0<r<\operatorname{dist}(x,\partial\Omega)$, relative to the agent at its center $x\in\Omega$:
\begin{align*}
F_r(x) := &
\int_{\partial B_r(x)} m(y) \left( -(v(y)-v(x)) \cdot \nu(y) \right)^{+} ds.
\end{align*}
The positive part retains only inward relative motion.

\begin{lemma}[Small-radius asymptotics of the influx]\label{lem:Fr_asymptotics}
Under Assumption~\ref{ass:local_regularity}, for every $x\in\Omega$, as $r\downarrow0$,
\begin{equation}\label{eq:Fr_asymptotics}
F_r(x)=r^2\,m(x)\int_0^{2\pi}\bigl(e(\theta)^\top Dv(x)\,e(\theta)\bigr)^{-}\,d\theta+O(r^{2+\alpha}),
\end{equation}
where, for every compact $K\Subset\Omega$, the $O(r^{2+\alpha})$ remainder is uniform for $x\in K$ and $0<r<\tfrac12\operatorname{dist}(K,\partial\Omega)$. Consequently,
\[
\lim_{r\to0}\frac{F_r(x)}{\pi r^2}=\frac{m(x)}{\pi}\int_0^{2\pi}\bigl(e(\theta)^\top Dv(x)\,e(\theta)\bigr)^{-}\,d\theta .
\]
\end{lemma}

\begin{proof}
Fix a compact set $K\Subset\Omega$ and set $\rho_K:=\tfrac12\operatorname{dist}(K,\partial\Omega)$, so that
\[
    K_{\rho_K}:=\{z\in\Omega:\operatorname{dist}(z,K)\le\rho_K\}\Subset\Omega .
\]
For $x\in K$ and $0<r<\rho_K$, $B_r(x)\subset K_{\rho_K}$, and with the parametrization of Section~\ref{sec:notation}, $\nu(y)=e(\theta)$ and $ds=r\,d\theta$. Since $v\in C^{1,\alpha}_{\mathrm{loc}}$, there is a constant $C_K$, depending only on the $C^{1,\alpha}$-seminorm of $v$ on $K_{\rho_K}$, with
\[
\bigl|v(y)-v(x)-Dv(x)(y-x)\bigr|\le C_K\,|y-x|^{1+\alpha}=C_K\,r^{1+\alpha}.
\]
Writing $a(\theta):=-r\,e(\theta)^\top Dv(x)\,e(\theta)$ and using $\nu(y)=e(\theta)$,
\[
-(v(y)-v(x))\cdot \nu(y)=a(\theta)+b(\theta),\qquad |b(\theta)|\le C_K\,r^{1+\alpha}.
\]
The map $s\mapsto s^{+}$ is $1$-Lipschitz, and $(-r\,t)^{+}=r\,t^{-}$ for $r>0$; hence
\[
\Bigl(-(v(y)-v(x))\cdot \nu(y)\Bigr)^{+}
= r\,\bigl(e(\theta)^\top Dv(x)\,e(\theta)\bigr)^{-}+O(r^{1+\alpha}),
\]
uniformly in $\theta$ and in $x\in K$. 
Since $m\in C^{0,1}_{\mathrm{loc}}$, we have $m(y)=m(x)+O(r)$,
with the remainder controlled by the Lipschitz seminorm of $m$ on
$K_{\rho_K}$. The quantity
$\bigl(e(\theta)^\top Dv(x)e(\theta)\bigr)^-$ is uniformly bounded
for $x\in K$ and $\theta\in[0,2\pi]$.
Multiplying and using $\alpha\le1$,
\[
m(y)\Bigl(-(v(y)-v(x))\cdot \nu(y)\Bigr)^{+}
= r\,m(x)\bigl(e(\theta)^\top Dv(x)\,e(\theta)\bigr)^{-}+O(r^{1+\alpha}).
\]
Integrating against $ds=r\,d\theta$, over $\theta\in[0,2\pi]$, yields \eqref{eq:Fr_asymptotics}, with a remainder constant depending only on $\rho_K$, $\|v\|_{C^{1,\alpha}(K_{\rho_K})}$, and $\|m\|_{C^{0,1}(K_{\rho_K})}$, hence uniform on $K$. Dividing by $\pi r^2$ and letting $r\to0$ gives the stated limit.
\end{proof}

We model local collision risk by combining the inward relative
flux per unit area, $\lim_{r\to0} F_r(x)/(\pi r^2)$, with
the local density of agents, $m(x)$.
Since
\[
    e(\theta)^\top Dv(x)e(\theta)
    =e(\theta)^\top A(x)e(\theta),
\]
where \(A(x):=\frac12(Dv(x)+Dv(x)^\top)\) is the
symmetric velocity gradient, we define the \emph{local collision intensity} by
\[
	\mathcal{C}(x) := m(x) \left( \lim_{r\to0}\frac{F_r(x)}{\pi r^2} \right) = \frac{m^2(x)}{\pi} \int_0^{2\pi} \left( e(\theta)^\top A(x) e(\theta) \right)^{-} d\theta.
\]
The factor $m(x)$ corresponds to the product approximation $m(x)m(y)$ for the density of ordered agent pairs, neglecting pair correlations, and the division by $\pi r^2$ normalizes out the interaction area, so $\mathcal C(x)$ is a normalized encounter intensity. Counting unordered pairs would introduce a factor $1/2$, which does not affect the minimization.
The local collision intensity is defined at every interior point
without prescribing boundary values. We define the total collision
functional, possibly taking the value \(+\infty\), by
\begin{equation}\label{eq:total_collision}
	\mathcal{C}_\Omega = \int_\Omega \mathcal{C}(x) dx
	= \frac{1}{\pi} \int_\Omega m^2(x) \int_0^{2\pi} \left( e(\theta)^\top A(x) e(\theta) \right)^{-} d\theta dx.
\end{equation}

We now compute the local collision intensity in terms of the eigenvalues of $A$. The small-radius inward-flux reduction to the symmetric velocity gradient follows the classical Saffman--Turner approximation \cite{saffmanCollisionDropsTurbulent1956}, and the resulting two-dimensional angular integral was evaluated in \cite[Appendix~A.1]{gustavssonCollisionsParticlesAdvected2008}; we include the short proof for completeness.

\begin{proposition}[Explicit form of the collision integral]\label{prop:collision_integral}
	Let $A$ be a real symmetric $2\times2$ matrix and set
\[
I(A)\;:=\;\int_{0}^{2\pi}\!\bigl(e(\theta)^{\top}A\,e(\theta)\bigr)^{-}\,d\theta.
\]
Let $\lambda_1,\lambda_2$ be the eigenvalues of $A$, labeled so that $\lambda_1\ge\lambda_2$.
Then $I(A)$ depends only on $(\lambda_1,\lambda_2)$ and
\[
I(A)=
\begin{cases}
0, & \lambda_1\ge\lambda_2\ge 0,\\[4pt]
-\pi\,(\lambda_1+\lambda_2), & 0\ge\lambda_1\ge\lambda_2,\\[6pt]
2\sqrt{-\lambda_1\lambda_2}\;-\;(\lambda_1+\lambda_2)\,
\arccos\!\Big(\dfrac{\lambda_1+\lambda_2}{\lambda_1-\lambda_2}\Big),
& \lambda_1>0>\lambda_2.
\end{cases}
\]
\end{proposition}
\begin{proof}
By the spectral theorem, $A=Q^{\top}\Lambda Q$ with $Q\in O(2)$ and $\Lambda=\diag(\lambda_1,\lambda_2)$. Then
\[
e(\theta)^{\top}(Q^{\top} \Lambda Q)\,e(\theta)=(Qe(\theta))^{\top} \Lambda \,(Qe(\theta)).
\]
In two dimensions, $Q$ is either a rotation or a reflection, so there exists $\phi\in\mathbb{R}$ such that $Qe(\theta)=e(\phi \pm \theta)$.
A change of variable $\vartheta=\phi \pm \theta$ and $2\pi$-periodicity yield
\[
I(A) = I(Q^{\top} \Lambda Q)=\int_{0}^{2\pi}\big(e(\phi \pm \theta)^{\top} \Lambda \,e(\phi \pm \theta)\big)^{-}d\theta
=\int_{0}^{2\pi}\big(e(\vartheta)^{\top} \Lambda \,e(\vartheta)\big)^{-}d\vartheta=I(\Lambda).
\]
Thus, $I(A)$ depends only on $(\lambda_1,\lambda_2)$.
In particular,
\[
I(A)=I(\Lambda)=\int_{0}^{2\pi}\big(\lambda_1\cos^2\theta+\lambda_2\sin^2\theta\big)^{-}\,d\theta.
\]

Next, set 
\[
f(\theta)=\lambda_1\cos^2\theta+\lambda_2\sin^2\theta
=\frac{\lambda_1+\lambda_2}{2}+\frac{\lambda_1-\lambda_2}{2}\cos(2\theta)
=:c+d\cos(2\theta).
\]
\begin{enumerate}
\item[i.] If $\lambda_1,\lambda_2\geq 0$, then $f(\theta)\geq 0$ for all $\theta$, hence $f(\theta)^-=0$ and $I(\Lambda)=0$.
\item[ii.] If $\lambda_1,\lambda_2\leq 0$, then $f(\theta)\leq0$ for all $\theta$, so $f(\theta)^-=-f(\theta)$ and
\[
I(\Lambda)
=-\int_{0}^{2\pi}\bigl(c+d\cos(2\theta)\bigr)\,d\theta
=-2\pi c
=-\pi(\lambda_1+\lambda_2).
\]
\item[iii.] If $\lambda_1>0>\lambda_2$, then $|c|<d$ (since $d^2-c^2=-\lambda_1\lambda_2>0$), so $f$ changes sign. Using $t=2\theta$ and the $2\pi$-periodicity of $\cos t$,
\[
I(\Lambda)=\int_0^{2\pi}\!\bigl(c+d\cos(2\theta)\bigr)^-\,d\theta
=\int_0^{2\pi}\!(c+d\cos t)^-\,dt.
\]
Let $\theta_0:=\arccos(-c/d)\in(0,\pi)$; then $\{t\in[0,2\pi]: c+d\cos t<0\}=(\theta_0,2\pi-\theta_0)$. Hence
\[
I(\Lambda)=\!\int_{\theta_0}^{2\pi-\theta_0}\!-(c+d\cos t)\,dt
=-2c(\pi-\theta_0)+2d\,\sin\theta_0.
\]
Since $\sin\theta_0=\sqrt{1-(c/d)^2}=\sqrt{d^2-c^2}/\,d$ and $\pi-\theta_0=\arccos(c/d)$, we obtain
\[
I(\Lambda)=2\sqrt{d^2-c^2}\;-\;2c\,\arccos\!\Big(\frac{c}{d}\Big)
=2\sqrt{-\lambda_1\lambda_2}\;-\;(\lambda_1+\lambda_2)\,
\arccos\!\Big(\frac{\lambda_1+\lambda_2}{\lambda_1-\lambda_2}\Big).
\]
\end{enumerate}
\end{proof}

\begin{remark}\label{rem:traceless_collision_integral}
In the traceless case $\tr A=0$, the eigenvalues satisfy $\lambda_1=-\lambda_2=\sqrt{-\det A}\ge0$, so the formula gives $I(A)=2\sqrt{-\det A}=2\lambda_1$. This is precisely the situation encountered in the dense crowd model of Section~\ref{sec:dense_flow}.
\end{remark}

\section{Dense incompressible stream-function formulation}
\label{sec:dense_flow}
%%%%%%%%%%%%%%%%%%%%%%%%%%%%%%%%%%%%%%%%%%%%%%%%%
%%%%%%%%%%%%%%%%%%%%%%%%%%%%%%%%%%%%%%%%%%%%%%%%%

We now specialize the collision functional to the dense incompressible regime
\[
    m\equiv 1,
    \qquad
    \div v=0
    \quad\text{in }\Omega .
\]

Here \(m\equiv1\) represents a fixed reference density in
normalized units. The variational analysis below is carried out
on bounded simply connected Lipschitz domains.

%%%%%%%%%%%%%%%%%%%%%%%%%%%%%%%%%%%%%%%%%%%%%%%%%
\subsection{Dense regime and boundary fluxes}
%%%%%%%%%%%%%%%%%%%%%%%%%%%%%%%%%%%%%%%%%%%%%%%%%

On a simply connected domain, every sufficiently smooth divergence-free
velocity field admits a stream function, unique up to an additive constant, such that
\[
    v=D^\perp\psi=(\psi_y,-\psi_x),
\]
and \(v\cdot\nu=\partial_\tau\psi\).

\begin{hyp}[Lipschitz domain and boundary flux]
\label{ass:dense_boundary}
The domain $\Omega\subset\mathbb R^2$ is bounded and simply
connected, with Lipschitz boundary $\Gamma=\partial\Omega$.
The prescribed signed boundary flux \(j\in L^1(\Gamma)\)
satisfies the compatibility condition
\begin{equation}\label{eq:dense_flux_compatibility}
    \int_\Gamma j\,ds=0 .
\end{equation}
\end{hyp}

The boundary condition \(v\cdot\nu=j\) prescribes inflow
on \(\{j<0\}\), outflow on \(\{j>0\}\), and zero normal
flow on \(\{j=0\}\), up to sets of arclength measure zero.

We parametrize \(\Gamma\) by arclength in the orientation induced by \(\tau\). If
\(j\in L^1(\Gamma)\) satisfies the compatibility condition
\eqref{eq:dense_flux_compatibility}, then there is a boundary primitive \(\bar j\in W^{1,1}(\Gamma)\), unique up to an additive constant, with \(\partial_\tau\bar j=j\) weakly on \(\Gamma\).
We fix the additive constant in \(\bar j\) once and for all.
 The boundary flux condition is then encoded by the Dirichlet condition in the trace sense:
\begin{equation}\label{eq:dense_stream_trace}
    \psi=\bar j     \quad\text{on }\Gamma .
\end{equation}
For nonsmooth stream functions, \eqref{eq:dense_stream_trace} is taken as the
meaning of the prescribed boundary flux.

%%%%%%%%%%%%%%%%%%%%%%%%%%%%%%%%%%%%%%%%%%%%%%%%%
\subsection{The dense collision energy}
%%%%%%%%%%%%%%%%%%%%%%%%%%%%%%%%%%%%%%%%%%%%%%%%%

Let \(\psi\in C^2(\overline\Omega)\). Since
\[
    D(D^\perp\psi)
    =
    \begin{pmatrix}
    \psi_{xy} & \psi_{yy}\\[1mm]
    -\psi_{xx} & -\psi_{xy}
    \end{pmatrix},
\]
the symmetric velocity gradient associated with \(D^\perp\psi\) is
\[
    A[\psi]
    :=
    \frac12\left(D(D^\perp\psi)+D(D^\perp\psi)^\top\right)
    =
    \begin{pmatrix}
    \psi_{xy} & \frac12(\psi_{yy}-\psi_{xx})\\[2mm]
    \frac12(\psi_{yy}-\psi_{xx}) & -\psi_{xy}
    \end{pmatrix}.
\]
This matrix is traceless, with eigenvalues
\[
    \lambda_\pm
    =
    \pm
    \lambda,
    \qquad
    \lambda
    :=
    \sqrt{
    \psi_{xy}^2+\frac14(\psi_{xx}-\psi_{yy})^2
    } .
\]
By Proposition~\ref{prop:collision_integral} and Remark~\ref{rem:traceless_collision_integral},
\[
    \int_0^{2\pi}
    \bigl(e(\theta)^\top A[\psi]e(\theta)\bigr)^-
    \,d\theta
    =
    2\lambda
    =
    \sqrt{
    4\psi_{xy}^2+(\psi_{xx}-\psi_{yy})^2
    }.
\]
Thus, in the dense incompressible regime, the total collision functional
\eqref{eq:total_collision} reduces to
\begin{equation}\label{eq:dense_energy_original}
    J[\psi]
    :=
    \frac1\pi
    \int_\Omega
    \sqrt{
    4\psi_{xy}^2+(\psi_{xx}-\psi_{yy})^2
    }\,dx\,dy .
\end{equation}

We introduce the second-order operator
\begin{equation}\label{eq:TD2_definition_dense}
    TD^2\psi
    :=
    \bigl(2\psi_{xy},\,\psi_{xx}-\psi_{yy}\bigr),
\end{equation}
so that
\begin{equation}\label{eq:dense_energy_TD2}
    J[\psi]
    =
    \frac1\pi
    \int_\Omega |TD^2\psi|\, dx\,dy ,
\end{equation}
a convex, positively one-homogeneous functional of \(TD^2\psi\).
Equivalently,
\[
    J[\psi]
    =\frac{\sqrt2}{\pi}
      \int_\Omega |E(D^\perp\psi)|\,dx\,dy.
\]
The integrand is affine along positive rays, so convexity alone does not give uniqueness; Proposition~\ref{prop:nonuniqueness} gives an explicit example with multiple minimizers.

The functional \(J\) controls only \(TD^2\psi\). In general, \(J\) does not control \(\|D^2\psi\|_{L^1(\Omega)}\), even for \(\psi\in C_c^\infty(\Omega)\); see \cite{ornsteinNoninequalityDifferentialOperators1962}. This distinction is essential
for the compactness and relaxation theory developed below.

%%%%%%%%%%%%%%%%%%%%%%%%%%%%%%%%%%%%%%%%%%%%%%%%%
\subsection{The admissible set}
%%%%%%%%%%%%%%%%%%%%%%%%%%%%%%%%%%%%%%%%%%%%%%%%%
\begin{hyp}[Finite-energy lifting]\label{ass:dense_lifting}
The boundary primitive \(\bar j\) associated with the data in
Assumption~\ref{ass:dense_boundary} admits a lifting satisfying
\begin{equation}\label{eq:dense_lifting_assumption}
    \psi_b\in W^{1,1}(\Omega),
    \qquad
    TD^2\psi_b\in L^1(\Omega;\mathbb R^2),
    \qquad
    \psi_b=\bar j
    \quad\text{on }\Gamma.
\end{equation}
Here \(TD^2\psi_b\) is understood in the sense of distributions
and the boundary equality is understood in the trace sense.
We fix one such lifting throughout the dense variational analysis.
\end{hyp}

For example, Assumption~\ref{ass:dense_lifting} holds whenever
\(\bar j\) is the boundary restriction of a \(C^2\) function
defined on a neighborhood of \(\overline\Omega\).
The admissible set is
\begin{equation}\label{eq:unrelaxed_class_A}
    \mathcal A
    :=
    \left\{
    \psi\in W^{1,1}(\Omega):
    TD^2\psi\in L^1(\Omega;\mathbb R^2),
    \quad
    \psi-\psi_b\in W^{1,1}_0(\Omega)
    \right\}.
\end{equation}

The class \(\mathcal A\) requires only the distributional combination \(TD^2\psi\) to be integrable; the relaxation will allow \(TD^2\psi\) to be a measure.
The dense variational problem is
\begin{equation}\label{eq:unrelaxed_dense_problem}
    \inf_{\psi\in\mathcal A} J[\psi].
\end{equation}
Since \(\psi_b\in\mathcal A\), \(0\le\inf_{\psi\in\mathcal A}J[\psi]\le J[\psi_b]<+\infty\).

%%%%%%%%%%%%%%%%%%%%%%%%%%%%%%%%%%%%%%%%%%%%%%%%%
\subsection{Geometry of the zero-trace kernel}
%%%%%%%%%%%%%%%%%%%%%%%%%%%%%%%%%%%%%%%%%%%%%%%%%

We characterize the zero-trace kernel of \(TD^2\).
It is trivial on domains other than disks. On a disk,
it is one-dimensional and corresponds to rigid rotation
about the center.

\begin{lemma}[Distributional kernel of \(TD^2\)]
\label{lem:dense_distributional_kernel}
Let \(\Omega\) be connected. If \(u\in W^{1,1}(\Omega)\)
satisfies \(TD^2u=0\) in distributions, then
\begin{equation}\label{eq:dense_distributional_kernel}
    u(x,y)
    =\frac{\alpha}{2}(x^2+y^2)+\beta x+\gamma y+\delta
\end{equation}
almost everywhere in \(\Omega\), for constants
\(\alpha,\beta,\gamma,\delta\).
\end{lemma}
\begin{proof}
All derivatives below are understood in the sense of distributions.
Differentiating \(u_{xy}=0\) gives
\[
    u_{xxy}=u_{xyy}=0.
\]
Differentiating \(u_{xx}-u_{yy}=0\) with respect to \(x\) and \(y\)
then yields
\[
    u_{xxx}=u_{xyy}=0,
    \qquad
    u_{yyy}=u_{xxy}=0.
\]
Thus both \(u_{xx}\) and \(u_{yy}\) have zero distributional
gradient. Since \(\Omega\) is connected and \(u_{xx}=u_{yy}\),
there is a constant \(\alpha\in\mathbb R\) such that
\[
    u_{xx}=u_{yy}=\alpha,
    \qquad
    u_{xy}=0.
\]

It follows that
\[
    D(u_x-\alpha x)=0,
    \qquad
    D(u_y-\alpha y)=0.
\]
Connectedness therefore gives constants \(\beta,\gamma\in\mathbb R\)
such that
\[
    u_x=\alpha x+\beta,
    \qquad
    u_y=\alpha y+\gamma.
\]
Consequently,
\[
    D\left(
        u-\frac{\alpha}{2}(x^2+y^2)-\beta x-\gamma y
      \right)=0.
\]
Applying the same constant-distribution property once more, we obtain
a constant \(\delta\in\mathbb R\) such that
\[
    u(x,y)
    =\frac{\alpha}{2}(x^2+y^2)+\beta x+\gamma y+\delta
    \quad\text{almost everywhere in }\Omega.
\]
This proves \eqref{eq:dense_distributional_kernel}.
\end{proof}

We define
\begin{equation}\label{eq:dense_zero_trace_kernel}
    \mathcal K_0
    :=\{q\in W^{1,1}_0(\Omega):TD^2q=0\}.
\end{equation}
Lemma~\ref{lem:dense_distributional_kernel} shows that
\(\mathcal K_0\) is finite-dimensional and consists of
polynomials of the form \eqref{eq:dense_distributional_kernel}.

\begin{proposition}[Geometry of the zero-trace kernel]
\label{prop:dense_kernel_geometry}
Let \(\Omega\subset\Rr^2\) be a bounded connected Lipschitz
domain. Then \(\mathcal K_0=\{0\}\) unless \(\Omega\) is an
open disk. If
\[
    \Omega=\{(x,y)\in\Rr^2:
              (x-x_0)^2+(y-y_0)^2<R^2\},
    \qquad R>0,
\]
then
\begin{equation}\label{eq:dense_disk_kernel_generator}
    \mathcal K_0=\operatorname{span}\{q_0\},
    \qquad
    q_0(x,y)=\frac{(x-x_0)^2+(y-y_0)^2-R^2}{2}.
\end{equation}
\end{proposition}
\begin{proof}
Suppose \(0\ne q\in\mathcal K_0\).
By Lemma~\ref{lem:dense_distributional_kernel},
\(q\) has the form \eqref{eq:dense_distributional_kernel}
and, by continuity and the zero trace, vanishes everywhere
on \(\Gamma\).
Its quadratic coefficient \(\alpha\) is nonzero, since an
affine function with zero boundary values on a bounded domain
must vanish identically.

Completing the square gives
\[
    q(x,y)=\frac{\alpha}{2}
    \bigl((x-x_0)^2+(y-y_0)^2+c\bigr)
\]
for some \((x_0,y_0)\in\Rr^2\) and \(c\in\Rr\).
Since \(q=0\) on \(\Gamma\) and \(\Gamma\) contains more than one point,
we must have \(c=-R^2\) for some \(R>0\), and hence
\(\Gamma\subset\partial B_R(x_0,y_0)\).

Let \(B=B_R(x_0,y_0)\), so that \(\Gamma=\partial\Omega\subset\partial B\).
For \(q_0\) defined by \eqref{eq:dense_disk_kernel_generator},
we have \(q=\alpha q_0\), so \(q_0=0\) on \(\partial\Omega\).
Since \(\Delta q_0=2>0\), the strong maximum principle gives
\(q_0<0\) throughout \(\Omega\), and hence \(\Omega\subset B\).

To see the reverse inclusion, note that \(B\cap\partial\Omega=\varnothing\)
because \(\partial\Omega\subset\partial B\). If there existed a point
\(z\in B\setminus\Omega\), then for any fixed \(p\in\Omega\subset B\),
the line segment connecting \(p\) to \(z\) would lie entirely in the
convex set \(B\) and would have to intersect \(\partial\Omega\),
contradicting \(B\cap\partial\Omega=\varnothing\). Thus \(B\subset\Omega\),
and consequently \(\Omega=B\).

Conversely, on this disk, the function \(q_0\) in
\eqref{eq:dense_disk_kernel_generator} belongs to \(\mathcal K_0\).
For any \(q\in\mathcal K_0\), subtracting \(\alpha q_0\),
where \(\alpha\) is its quadratic coefficient, leaves an affine
function with zero boundary values, hence the zero function.
Therefore \(\mathcal K_0=\operatorname{span}\{q_0\}\).
\end{proof}

%%%%%%%%%%%%%%%%%%%%%%%%%%%%%%%%%%%%%%%%%%%%%%%%%
\subsection{Normalization and coercivity}
\label{sec:dense_coercivity}
%%%%%%%%%%%%%%%%%%%%%%%%%%%%%%%%%%%%%%%%%%%%%%%%%

Coercivity holds after removing the component in the
zero-trace kernel \(\mathcal K_0\).
The compactness argument uses the relation between
\(TD^2\psi\) and the symmetric gradient of \(D^\perp\psi\).
For every \(\psi\in W^{1,1}(\Omega)\) with
\(TD^2\psi\in\mathcal M(\Omega;\Rr^2)\), a direct
computation gives \(D^\perp\psi\in BD(\Omega)\) and
\begin{equation}\label{eq:dense_BD_norm}
    \|D^\perp\psi\|_{BD(\Omega)}
    =
    \|D\psi\|_{L^1(\Omega)}
    +\frac1{\sqrt2}|TD^2\psi|(\Omega).
\end{equation}
We use the compact embedding
\(BD(\Omega)\hookrightarrow\hookrightarrow
L^1(\Omega;\Rr^2)\) on bounded Lipschitz domains;
see \cite[Theorem~2.10]{barrosoRelaxationTheoremSpace2000}.

\begin{example}[Failure of coercivity on a disk]
\label{ex:dense_disk_kernel}
Let \(\Omega\) be a disk and let \(q_0\) be given by
\eqref{eq:dense_disk_kernel_generator}.
For the fixed lifting \(\psi_b\) in
Assumption~\ref{ass:dense_lifting}, the sequence
\(\psi_n=\psi_b+nq_0\) belongs to \(\mathcal A\) and
satisfies
\[
    J[\psi_n]=J[\psi_b],
    \qquad
    \|\psi_n\|_{W^{1,1}(\Omega)}\longrightarrow+\infty.
\]
Thus the collision energy does not control the
\(W^{1,1}(\Omega)\) norm, even when the boundary trace
is fixed.
\end{example}

\begin{proposition}[Coercivity]
\label{prop:dense_coercivity}
Let \(\Omega\subset\Rr^2\) be a bounded connected Lipschitz
domain. Define
\begin{equation}\label{eq:dense_kernel_projection}
    P_0u:=
    \begin{cases}
        \displaystyle
        \frac{\int_\Omega uq_0\,dx\,dy}
             {\int_\Omega q_0^2\,dx\,dy}\,q_0,
        & \text{if \(\Omega\) is a disk},\\[2mm]
        0,
        & \text{otherwise},
    \end{cases}
\end{equation}
where \(q_0\) is given by
\eqref{eq:dense_disk_kernel_generator}.
Then \(P_0:W^{1,1}(\Omega)\to\mathcal K_0\) is a
continuous linear projection, and there is a constant
\(C_0>0\), depending only on \(\Omega\), such that
\begin{equation}\label{eq:dense_zero_trace_coercivity}
    \|u-P_0u\|_{W^{1,1}(\Omega)}
    \leq C_0|TD^2u|(\Omega)
\end{equation}
for every \(u\in W^{1,1}_0(\Omega)\) with
\(TD^2u\in\mathcal M(\Omega;\Rr^2)\).

If \(\psi_b\in W^{1,1}(\Omega)\) with
\(TD^2\psi_b\in L^1(\Omega;\Rr^2)\) is fixed, then
there is a constant \(C>0\), depending only on
\(\Omega\) and \(\psi_b\), such that
\begin{equation}\label{eq:dense_coercivity}
    \|\psi-P_0\psi\|_{W^{1,1}(\Omega)}
    \leq C\bigl(1+|TD^2\psi|(\Omega)\bigr)
\end{equation}
for every \(\psi\in W^{1,1}(\Omega)\) satisfying
\(TD^2\psi\in\mathcal M(\Omega;\Rr^2)\) and
\(\psi-\psi_b\in W^{1,1}_0(\Omega)\).
\end{proposition}

\begin{proof}
By Proposition~\ref{prop:dense_kernel_geometry},
\(P_0\) is a projection onto \(\mathcal K_0\).
Its continuity is immediate when \(P_0=0\), and follows
from boundedness of \(q_0\) in the disk case.

Suppose \eqref{eq:dense_zero_trace_coercivity} fails.
Subtracting the projection and normalizing gives a sequence
\(z_n\in W^{1,1}_0(\Omega)\) such that
\[
    \|z_n\|_{W^{1,1}(\Omega)}=1,
    \qquad
    P_0z_n=0,
    \qquad
    |TD^2z_n|(\Omega)\longrightarrow0.
\]
By \eqref{eq:dense_BD_norm}, the fields
\(D^\perp z_n\) are bounded in \(BD(\Omega)\).
The compact embedding therefore gives a subsequence
for which \(Dz_n\) converges strongly in \(L^1\).
Poincar\'e's inequality then gives strong convergence
in \(W^{1,1}(\Omega)\) to some
\(z\in W^{1,1}_0(\Omega)\).
Distributional convergence and continuity of \(P_0\) yield
\[
    \|z\|_{W^{1,1}(\Omega)}=1,
    \qquad
    TD^2z=0,
    \qquad
    P_0z=0.
\]
Thus \(z\in\mathcal K_0\), so \(z=P_0z=0\),
a contradiction.

For the prescribed-trace estimate, apply
\eqref{eq:dense_zero_trace_coercivity} to
\(u=\psi-\psi_b\). By linearity of \(P_0\) and
the triangle inequality,
\[
    \|\psi-P_0\psi\|_{W^{1,1}(\Omega)}
    \leq
    \|\psi_b-P_0\psi_b\|_{W^{1,1}(\Omega)}
    +C_0\bigl(
        |TD^2\psi|(\Omega)
        +\|TD^2\psi_b\|_{L^1(\Omega)}
      \bigr).
\]
Absorbing the fixed lifting terms into \(C\) proves
\eqref{eq:dense_coercivity}.
\end{proof}

For \(\psi\) satisfying the prescribed-trace hypotheses
of Proposition~\ref{prop:dense_coercivity}, define
\[
    \widehat\psi:=\psi-P_0\psi.
\]
Since \(P_0\psi\in\mathcal K_0\),
\begin{equation}\label{eq:dense_normalization_invariance}
    \widehat\psi-\psi_b\in W^{1,1}_0(\Omega),
    \qquad
    TD^2\widehat\psi=TD^2\psi,
    \qquad
    P_0\widehat\psi=0.
\end{equation}
In particular, if \(\psi\in\mathcal A\), then
\(\widehat\psi\in\mathcal A\) and
\(J[\widehat\psi]=J[\psi]\).
When \(\Omega\) is not a disk, the normalization is
the identity.

\begin{remark}[Interpretation of the disk normalization]
\label{rem:dense_normalization_interpretation}
On a disk, the normalized stream function satisfies
\[
    \int_\Omega\widehat\psi q_0\,dx\,dy=0.
\]
Among the functions in \(\psi+\mathcal K_0\),
\(\widehat\psi\) is the unique one satisfying this condition.
Since
\[
    D^\perp q_0=(y-y_0,-(x-x_0)),
\]
the normalization subtracts a rigid rotation tangent to
the circular boundary. It may change the velocity, but
preserves the prescribed trace and the entire measure
\(TD^2\psi\). It is used to obtain compactness and does
not select a preferred physical flow.
Only zero-trace kernel elements are subtracted; a general
element of \(\ker TD^2\) need not preserve the prescribed
trace or boundary flux.
\end{remark}

%%%%%%%%%%%%%%%%%%%%%%%%%%%%%%%%%%%%%%%%%%%%%%%%%
%%%%%%%%%%%%%%%%%%%%%%%%%%%%%%%%%%%%%%%%%%%%%%%%%
\section{The relaxed dense variational problem}
\label{sec:relaxed_dense_problem}
%%%%%%%%%%%%%%%%%%%%%%%%%%%%%%%%%%%%%%%%%%%%%%%%%
%%%%%%%%%%%%%%%%%%%%%%%%%%%%%%%%%%%%%%%%%%%%%%%%%

Proposition~\ref{prop:dense_coercivity} bounds every
minimizing sequence for \(J\) in \(W^{1,1}(\Omega)\)
after the normalization \(\psi\mapsto\psi-P_0\psi\).
This normalization preserves the prescribed trace and
energy and is the identity when \(\Omega\) is not a disk.

The energy bounds the \(L^1\) norm of \(TD^2\psi\), but
these derivatives may concentrate along sequences of
bounded energy. We therefore retain the prescribed trace
and allow \(TD^2\psi\) to be a finite vector-valued Radon
measure. Compactness and existence follow from the normalized
coercivity estimate. The recovery construction and
the identification of the relaxation require no
normalization.

%%%%%%%%%%%%%%%%%%%%%%%%%%%%%%%%%%%%%%%%%%%%%%%%%
\subsection{The relaxed admissible class and energy}
%%%%%%%%%%%%%%%%%%%%%%%%%%%%%%%%%%%%%%%%%%%%%%%%%

We define
\begin{equation}\label{eq:XT_definition}
    \mathcal X_T(\Omega)
    :=
    \left\{
    \psi\in W^{1,1}(\Omega):
    TD^2\psi\in\mathcal M(\Omega;\mathbb R^2)
    \right\}.
\end{equation}
The operator \(TD^2\) is understood in the sense of
distributions. Thus, for
\(\Phi=(\Phi_1,\Phi_2)\in C_c^\infty(\Omega;\mathbb R^2)\),
\[
    \langle TD^2\psi,\Phi\rangle
    =
    \int_\Omega
    \psi\left(
        2\partial_{xy}\Phi_1
        +(\partial_{xx}-\partial_{yy})\Phi_2
    \right)\,dx\,dy.
\]

The boundary flux is imposed through the same trace condition
as in the unrelaxed problem. We define
\begin{equation}\label{eq:relaxed_class_B}
    \mathcal B
    :=
    \left\{
    \psi\in\mathcal X_T(\Omega):
    \psi-\psi_b\in W^{1,1}_0(\Omega)
    \right\}.
\end{equation}
The boundary condition is equivalent to
\(\psi=\bar j\) almost everywhere on
\(\Gamma\) in the trace sense. By Assumption~\ref{ass:dense_lifting},
\(\psi_b\in\mathcal A\subset\mathcal B\).

We extend the functional in \eqref{eq:dense_energy_TD2}
to \(\mathcal X_T(\Omega)\) by setting
\[
    \mathcal J[\psi]
    :=
    \frac1\pi |TD^2\psi|(\Omega).
\]
If
\[
    TD^2\psi=g_\psi\,\mathcal L^2+\mu_\psi^s
\]
is the Lebesgue decomposition of \(TD^2\psi\), then
\[
    \mathcal J[\psi]
    =
    \frac1\pi\int_\Omega |g_\psi|\,dx\,dy
    +\frac1\pi|\mu_\psi^s|(\Omega).
\]
For a relaxed stream function, the natural collision measure
is \(\pi^{-1}|TD^2\psi|\). Its absolutely continuous density
is \(\pi^{-1}|g_\psi|\), while its singular part records
concentration in the relaxed limit.
Thus \(\mathcal J\) extends the original dense energy:
if \(\psi\in\mathcal A\), then \(TD^2\psi\) is absolutely
continuous and
\[
    \mathcal J[\psi]
    =
    \frac1\pi\int_\Omega |TD^2\psi|\,dx\,dy
    =
    J[\psi].
\]

The following lemma records the behavior of the energy
along a sequence converging to a prescribed limit.
It applies on every domain satisfying the standing
assumptions and does not use normalization.

\begin{lemma}[Limits of sequences with bounded energy]
\label{lem:relaxed_dense_closedness}
Under Assumptions~\ref{ass:dense_boundary}
and~\ref{ass:dense_lifting}, let \(\{\psi_n\}\subset\mathcal B\) satisfy
\[
    \psi_n\to\psi\quad\text{in }W^{1,1}(\Omega),
    \qquad
    \sup_n\mathcal J[\psi_n]<+\infty.
\]
Then \(\psi\in\mathcal B\),
\[
    TD^2\psi_n\cws TD^2\psi
    \quad\text{in }\mathcal M(\Omega;\Rr^2),
    \qquad
    \mathcal J[\psi]
    \leq\liminf_{n\to\infty}\mathcal J[\psi_n].
\]
Moreover, the weak-\(*\) convergence holds for the entire sequence.
\end{lemma}

\begin{proof}
The energy bound gives a weak-\(*\) convergent subsequence
of \(TD^2\psi_n\). Strong \(W^{1,1}\) convergence implies
distributional convergence, so the limiting measure is
\(TD^2\psi\). Thus \(\psi\in\mathcal X_T(\Omega)\).
Since \(W^{1,1}_0(\Omega)\) is closed and
\(\psi_n-\psi_b\in W^{1,1}_0(\Omega)\), we also have
\(\psi-\psi_b\in W^{1,1}_0(\Omega)\). Hence
\(\psi\in\mathcal B\).

Every subsequence of \(TD^2\psi_n\) has a further
weak-\(*\) convergent subsequence by the uniform energy
bound. Distributional convergence identifies every such
limit with \(TD^2\psi\). 
Hence the entire sequence converges weak-\(*\) to \(TD^2\psi\).
For every \(\Phi\in C_c(\Omega;\mathbb R^2)\) with
\(\|\Phi\|_\infty\le1\), we have
\[
    \int_\Omega \Phi\cdot d(TD^2\psi)
    =
    \lim_{n\to\infty}
    \int_\Omega \Phi\cdot d(TD^2\psi_n)
    \le
    \liminf_{n\to\infty}|TD^2\psi_n|(\Omega).
\]
Taking the supremum over such \(\Phi\), by the dual
characterization of the total variation of a vector-valued
Radon measure, gives
\[
    |TD^2\psi|(\Omega)
    \le
    \liminf_{n\to\infty}|TD^2\psi_n|(\Omega).
\]
Dividing by \(\pi\) proves the energy inequality.
\end{proof}

Theorem~\ref{thm:relaxed_dense_relaxation_identity}
will identify \(\mathcal J\) on \(\mathcal B\) with
the relaxation of \(J\) on \(\mathcal A\) for strong
\(W^{1,1}(\Omega)\) convergence.

%%%%%%%%%%%%%%%%%%%%%%%%%%%%%%%%%%%%%%%%%%%%%%%%%
\subsection{Compactness and existence}
%%%%%%%%%%%%%%%%%%%%%%%%%%%%%%%%%%%%%%%%%%%%%%%%%

We prove compactness after normalization using
Proposition~\ref{prop:dense_coercivity} and the compact
embedding of \(BD(\Omega)\) into \(L^1(\Omega;\Rr^2)\).
The direct method then gives a minimizer of \(\mathcal J\)
over the full class \(\mathcal B\).

\begin{proposition}[Compactness after normalization]
\label{prop:relaxed_dense_compactness}
Under Assumptions~\ref{ass:dense_boundary}
and~\ref{ass:dense_lifting}, let
\(\{\psi_n\}\subset\mathcal B\) satisfy
\[
    \sup_n|TD^2\psi_n|(\Omega)<+\infty.
\]
Set
\[
    \widehat\psi_n:=\psi_n-P_0\psi_n,
\]
where \(P_0\) is given by
\eqref{eq:dense_kernel_projection}.
Then, for every \(n\),
\begin{equation}\label{eq:relaxed_normalization_invariance}
    \widehat\psi_n\in\mathcal B,
    \qquad
    TD^2\widehat\psi_n=TD^2\psi_n,
    \qquad
    \mathcal J[\widehat\psi_n]=\mathcal J[\psi_n].
\end{equation}
There are \(\psi\in\mathcal B\), with \(P_0\psi=0\),
and a subsequence, not relabeled, such that
\[
    \widehat\psi_n\to\psi
    \quad\text{in }W^{1,1}(\Omega),
    \qquad
    TD^2\psi_n\cws TD^2\psi
    \quad\text{in }\mathcal M(\Omega;\Rr^2).
\]
If \(\Omega\) is not a disk, then
\(\widehat\psi_n=\psi_n\).
\end{proposition}

\begin{proof}
Equation~\eqref{eq:dense_normalization_invariance} gives
\eqref{eq:relaxed_normalization_invariance} and
\(P_0\widehat\psi_n=0\).
Proposition~\ref{prop:dense_coercivity} bounds
\(\widehat\psi_n\) in \(W^{1,1}(\Omega)\).
Together with \eqref{eq:dense_BD_norm} and the energy
bound, this shows that \(D^\perp\widehat\psi_n\) is
bounded in \(BD(\Omega)\).
The compact embedding gives a subsequence for which
\(D\widehat\psi_n\) converges strongly in \(L^1\).

Since \(\widehat\psi_n-\widehat\psi_m\in W^{1,1}_0(\Omega)\),
Poincar\'e's inequality gives
\[
    \|\widehat\psi_n-\widehat\psi_m\|_{L^1(\Omega)}
    \leq C
    \|D\widehat\psi_n-D\widehat\psi_m\|_{L^1(\Omega)}.
\]
Thus this subsequence converges strongly in
\(W^{1,1}(\Omega)\) to some \(\psi\).
Lemma~\ref{lem:relaxed_dense_closedness} gives
\(\psi\in\mathcal B\) and the asserted weak-\(*\)
convergence, since \(TD^2\widehat\psi_n=TD^2\psi_n\).
Finally, continuity of \(P_0\) gives \(P_0\psi=0\).
\end{proof}

\begin{corollary}[Compactness of relaxed energy sublevels]
\label{cor:relaxed_sublevels}
Under Assumptions~\ref{ass:dense_boundary}
and~\ref{ass:dense_lifting}, for \(M\ge0\), define
\[
    \mathcal B_M
    :=\{\psi\in\mathcal B:\mathcal J[\psi]\le M\},
    \qquad
    \widehat{\mathcal B}_M
    :=\{\psi\in\mathcal B_M:P_0\psi=0\}.
\]
Then \(\widehat{\mathcal B}_M\) is compact in the strong
\(W^{1,1}(\Omega)\) topology.
If \(\Omega\) is not a disk, then
\(\widehat{\mathcal B}_M=\mathcal B_M\).
If \(\Omega\) is a disk, every nonempty
\(\mathcal B_M\) is unbounded in \(W^{1,1}(\Omega)\).
\end{corollary}

\begin{proof}
Let \(\{\psi_n\}\subset\widehat{\mathcal B}_M\).
Since \(P_0\psi_n=0\) for every \(n\), the normalized
sequence satisfies \(\widehat\psi_n=\psi_n\).
Proposition~\ref{prop:relaxed_dense_compactness} gives
a strongly convergent subsequence whose limit belongs
to \(\mathcal B\) and satisfies \(P_0\psi=0\).
Lemma~\ref{lem:relaxed_dense_closedness} preserves the
energy bound, proving compactness of
\(\widehat{\mathcal B}_M\).
If \(\Omega\) is not a disk, then \(P_0=0\), so
\(\widehat{\mathcal B}_M=\mathcal B_M\).
If \(\Omega\) is a disk and \(\psi\in\mathcal B_M\), then
\(\psi+tq_0\in\mathcal B_M\) for every \(t\in\Rr\),
because \(q_0\) has zero trace and \(TD^2q_0=0\).
However,
\[
    \|\psi+tq_0\|_{W^{1,1}(\Omega)}
    \ge |t|\,\|q_0\|_{W^{1,1}(\Omega)}
          -\|\psi\|_{W^{1,1}(\Omega)}
    \longrightarrow+\infty
    \quad\text{as }|t|\to\infty.
\]
Thus every nonempty unnormalized sublevel on a disk is
unbounded.
\end{proof}

\begin{theorem}[Existence of relaxed minimizers]
\label{thm:relaxed_dense_existence}
Under Assumptions~\ref{ass:dense_boundary}
and~\ref{ass:dense_lifting}, \(\mathcal J\) attains
its minimum on \(\mathcal B\).
\end{theorem}

\begin{proof}
Since \(\psi_b\in\mathcal B\) and
\(\mathcal J[\psi_b]=J[\psi_b]<+\infty\), the infimum
is finite and nonnegative.
Let \(\{\psi_n\}\subset\mathcal B\) be a minimizing
sequence. Then \(|TD^2\psi_n|(\Omega)=\pi\mathcal J[\psi_n]\)
is uniformly bounded.

By Proposition~\ref{prop:relaxed_dense_compactness},
the normalized sequence
\(\widehat\psi_n=\psi_n-P_0\psi_n\) has a subsequence
converging strongly in \(W^{1,1}(\Omega)\) to some
\(\psi\in\mathcal B\).
Since normalization preserves the energy,
Lemma~\ref{lem:relaxed_dense_closedness} gives
\[
    \mathcal J[\psi]
    \leq\liminf_{n\to\infty}\mathcal J[\widehat\psi_n]
    =\liminf_{n\to\infty}\mathcal J[\psi_n]
    =\inf_{\varphi\in\mathcal B}\mathcal J[\varphi].
\]
Since \(\psi\in\mathcal B\), it is a minimizer.
\end{proof}

%%%%%%%%%%%%%%%%%%%%%%%%%%%%%%%%%%%%%%%%%%%%%%%%%
\subsection{Strict recovery with the prescribed trace}
%%%%%%%%%%%%%%%%%%%%%%%%%%%%%%%%%%%%%%%%%%%%%%%%%

We approximate every \(\psi\in\mathcal B\) by functions in
\(\mathcal A\) with the same prescribed trace and convergent
collision energy. The construction applies interior
mollification directly to \(\psi\) and combines the local
approximations using a partition of unity. No normalization
is needed, whether or not \(\Omega\) is a disk.

\begin{theorem}[Boundary-preserving strict \(TD^2\)-recovery]
\label{thm:TD_strict_recovery}
Under Assumptions~\ref{ass:dense_boundary}
and~\ref{ass:dense_lifting}, for every \(\psi\in\mathcal B\),
there exists a sequence
\(\{\psi_k\}_{k\in\Nn}\subset\mathcal A\cap C^\infty(\Omega)\)
such that
\[
    \psi_k-\psi\in W^{1,1}_0(\Omega),
    \qquad
    \psi_k\to\psi
    \quad\text{in }W^{1,1}(\Omega),
\]
\[
    TD^2\psi_k\cws TD^2\psi
    \quad\text{in }\mathcal M(\Omega;\Rr^2),
\]
and
\begin{equation}\label{eq:TD_strict_recovery_variation}
    |TD^2\psi_k|(\Omega)\to|TD^2\psi|(\Omega).
\end{equation}
Consequently, \(J[\psi_k]\to\mathcal J[\psi]\).
\end{theorem}

\begin{proof}
\noindent\textbf{Local smoothing and the product rule.}
Set \(\mu:=TD^2\psi\).
Choose a nonnegative, locally finite partition of unity
\(\{\varphi_i\}_{i\ge1}\subset C_c^\infty(\Omega)\),
with \(\sum_i\varphi_i=1\) on \(\Omega\).
Extend each \(\varphi_i\) smoothly by zero to \(\Rr^2\).
Denote by \(\widetilde\psi\) and \(\widetilde\mu\) the
zero extensions of \(\psi\) as an \(L^1\) function and of
\(\mu\) as a finite measure.
The extension of \(\psi\) is used only for interior
convolution; no derivative identity across \(\Gamma\)
is required.

Let \(\rho\) be a nonnegative even standard mollifier
supported in \(B_1(0)\), and set
\[
    \psi_\delta:=\rho_\delta*\widetilde\psi.
\]
As \(\delta\downarrow0\), we have
\(\psi_\delta\to\psi\) in
\(W^{1,1}_{\mathrm{loc}}(\Omega)\).
For
\[
    0<\delta<
    \frac12\operatorname{dist}
    (\supp\varphi_i,\Rr^2\setminus\Omega),
\]
the identity
\begin{equation}\label{eq:dense_local_convolution}
    TD^2\psi_\delta=\rho_\delta*\widetilde\mu
\end{equation}
holds on a neighborhood of \(\supp\varphi_i\).

For \(\varphi\in C_c^\infty(\Omega)\) and
\(q\in W^{1,1}(\Omega)\), define
\begin{align*}
    \mathcal R_\varphi[q]:=\bigl(&
    2\varphi_xq_y+2\varphi_yq_x+2\varphi_{xy}q,\\
    &2\varphi_xq_x-2\varphi_yq_y
      +(\varphi_{xx}-\varphi_{yy})q\bigr).
\end{align*}
The distributional product rule gives
\begin{equation}\label{eq:dense_recovery_product}
    TD^2(\varphi q)
    =\varphi TD^2q+\mathcal R_\varphi[q].
\end{equation}
The remainder involves only \(q\) and its first derivatives.

For each \(k,i\ge1\), choose \(\delta_{k,i}>0\) small enough
for \eqref{eq:dense_local_convolution} to hold near
\(\supp\varphi_i\) and for all three inequalities
\begin{align}
    \|\varphi_i(\psi_{\delta_{k,i}}-\psi)\|_{W^{1,1}(\Omega)}
    &\leq \frac{2^{-i}}{k},
    \label{eq:dense_recovery_choice1}\\
    \|\mathcal R_{\varphi_i}
        [\psi_{\delta_{k,i}}-\psi]\|_{L^1(\Omega)}
    &\leq \frac{2^{-i}}{k},
    \label{eq:dense_recovery_choice2}\\
    \|\rho_{\delta_{k,i}}*\varphi_i-\varphi_i
        \|_{L^\infty(\Rr^2)}
    &\leq \frac{2^{-i}}{k}.
    \label{eq:dense_recovery_choice3}
\end{align}
The first two follow from local \(W^{1,1}\) convergence.
The third follows from uniform convergence of mollifications
of the fixed smooth compactly supported function \(\varphi_i\).

\medskip
\noindent\textbf{Preservation of the scalar trace.}
Define
\begin{equation}\label{eq:dense_interior_approximation}
    \psi_k
    :=\sum_{i=1}^\infty
       \varphi_i\psi_{\delta_{k,i}}.
\end{equation}
The sum is locally finite, so \(\psi_k\in C^\infty(\Omega)\).
Moreover,
\[
    \psi_k-\psi
    =\sum_i\varphi_i(\psi_{\delta_{k,i}}-\psi).
\]
By \eqref{eq:dense_recovery_choice1}, this series converges
absolutely in \(W^{1,1}(\Omega)\), with
\[
    \|\psi_k-\psi\|_{W^{1,1}(\Omega)}\leq\frac1k.
\]
Every summand belongs to \(W^{1,1}_0(\Omega)\), since it
belongs to \(W^{1,1}(\Omega)\) and has compact support
in \(\Omega\). Closedness of \(W^{1,1}_0(\Omega)\)
therefore gives
\[
    \psi_k-\psi\in W^{1,1}_0(\Omega),
    \qquad
    \psi_k\to\psi
    \quad\text{in }W^{1,1}(\Omega).
\]
Since \(\psi\in\mathcal B\), we also have
\(\psi_k-\psi_b\in W^{1,1}_0(\Omega)\).
Thus the prescribed scalar trace is preserved exactly.

\medskip
\noindent\textbf{Convergence of the energy.}
Local finiteness and \(\sum_i\varphi_i=1\) give
\[
    \sum_iD\varphi_i=0,
    \qquad
    \sum_iD^2\varphi_i=0.
\]
Hence \(\sum_i\mathcal R_{\varphi_i}[\psi]=0\) locally.
By local finiteness, we may differentiate the sum defining
\(\psi_k\) term by term. The product rule
\eqref{eq:dense_recovery_product} and the local convolution
identity \eqref{eq:dense_local_convolution} give
\begin{align*}
    TD^2\psi_k
    &=
    \sum_i TD^2\bigl(\varphi_i\psi_{\delta_{k,i}}\bigr)\\
    &=
    \sum_i\varphi_i
        \bigl(\rho_{\delta_{k,i}}*\widetilde\mu\bigr)
    +\sum_i\mathcal R_{\varphi_i}
        [\psi_{\delta_{k,i}}].
\end{align*}
Since \(\mathcal R_{\varphi_i}\) is linear in its argument
and \(\sum_i\mathcal R_{\varphi_i}[\psi]=0\), we have
\[
    \sum_i\mathcal R_{\varphi_i}[\psi_{\delta_{k,i}}]
    =
    \sum_i\mathcal R_{\varphi_i}
        [\psi_{\delta_{k,i}}-\psi].
\]
Thus, in distributions,
\[
    TD^2\psi_k=M_k+r_k,
\]
where
\[
    M_k
    :=\sum_i\varphi_i
       (\rho_{\delta_{k,i}}*\widetilde\mu),
    \qquad
    r_k
    :=\sum_i\mathcal R_{\varphi_i}
       [\psi_{\delta_{k,i}}-\psi].
\]
By \eqref{eq:dense_recovery_choice2}, the series defining
\(r_k\) converges absolutely in \(L^1(\Omega;\Rr^2)\), and
\[
    \|r_k\|_{L^1(\Omega)}\leq\frac1k.
\]

Since \(\rho\) is even,
\[
    \int_{\Rr^2}
        \varphi_i(w)\rho_{\delta_{k,i}}(w-z)\,dw
    =
    (\rho_{\delta_{k,i}}*\varphi_i)(z).
\]
Note also that 
\eqref{eq:dense_recovery_choice3} gives, for every
\(z\in\Omega\),
\[
    (\rho_{\delta_{k,i}}*\varphi_i)(z)
    \leq \varphi_i(z)+\frac{2^{-i}}k.
\]
Summing over \(i\) and using
\(\sum_i\varphi_i(z)=1\) and \(\sum_{i=1}^\infty2^{-i}=1\),
we obtain
\[
    \sum_i(\rho_{\delta_{k,i}}*\varphi_i)(z)
    \leq 1+\frac1k.
\]
Using the convolution identity, the preceding bound,
nonnegativity, and Tonelli's theorem, we obtain
\begin{align*}
    \int_\Omega|M_k|\,dx\,dy
    &\leq
    \sum_i\int_{\Rr^2}\varphi_i
        (\rho_{\delta_{k,i}}*|\widetilde\mu|)
        \,dx\,dy\\
    &=
    \int_\Omega
        \sum_i(\rho_{\delta_{k,i}}*\varphi_i)(z)
        \,d|\mu|(z)\\
    &\leq
    \left(1+\frac1k\right)|\mu|(\Omega).
\end{align*}

Consequently,
\begin{equation}\label{eq:dense_recovery_upper}
    \int_\Omega|TD^2\psi_k|\,dx\,dy
    \leq
    \left(1+\frac1k\right)|\mu|(\Omega)+\frac1k.
\end{equation}
In particular, \(TD^2\psi_k\in L^1(\Omega;\Rr^2)\).
Together with the trace condition, this gives
\(\psi_k\in\mathcal A\).

The bound \eqref{eq:dense_recovery_upper} and
Lemma~\ref{lem:relaxed_dense_closedness} give
\[
    TD^2\psi_k\cws TD^2\psi
    \quad\text{in }\mathcal M(\Omega;\Rr^2),
\]
and lower semicontinuity yields
\[
    |\mu|(\Omega)
    \leq
    \liminf_{k\to\infty}
    \int_\Omega|TD^2\psi_k|\,dx\,dy.
\]
Together with \eqref{eq:dense_recovery_upper}, this proves
\eqref{eq:TD_strict_recovery_variation}.
Finally,
\[
    J[\psi_k]
    =\frac1\pi|TD^2\psi_k|(\Omega)
    \longrightarrow
    \frac1\pi|TD^2\psi|(\Omega)
    =\mathcal J[\psi].
\]
\end{proof}

The approximation argument applies to every
\(\psi\in\mathcal X_T(\Omega)\) on an arbitrary bounded
open set \(\Omega\) and gives
\(\psi_k-\psi\in W^{1,1}_0(\Omega)\).
On a Lipschitz domain, this means that the scalar traces
are equal.

The approximants are smooth in \(\Omega\) and have
integrable \(TD^2\psi_k\). Neither integrability of the
full Hessian nor smoothness up to \(\partial\Omega\)
is asserted. This is precisely the regularity needed
for \(\mathcal A\).

%%%%%%%%%%%%%%%%%%%%%%%%%%%%%%%%%%%%%%%%%%%%%%%%%
\subsection{Identification of the relaxation}
%%%%%%%%%%%%%%%%%%%%%%%%%%%%%%%%%%%%%%%%%%%%%%%%%

For \(\psi\in W^{1,1}(\Omega)\), define
\begin{equation}\label{eq:J_rel_definition}
    J_{\rm rel}[\psi]
    :=
    \inf
    \left\{
    \liminf_{k\to\infty}J[\psi_k]:
    \psi_k\in\mathcal A,
    \quad
    \psi_k\to\psi
    \text{ in }W^{1,1}(\Omega)
    \right\}.
\end{equation}
We use the convention \(\inf\varnothing=+\infty\).
We also set
\[
    \overline{\mathcal J}[\psi]
    :=
    \begin{cases}
    \mathcal J[\psi], & \psi\in\mathcal B,\\[1mm]
    +\infty, & \psi\notin\mathcal B.
    \end{cases}
\]
The following theorem identifies
\(\overline{\mathcal J}\) with the relaxation of \(J\)
on \(\mathcal A\) in the strong \(W^{1,1}(\Omega)\)
topology.

\begin{theorem}[Relaxation identity]
\label{thm:relaxed_dense_relaxation_identity}
Under Assumptions~\ref{ass:dense_boundary}
and~\ref{ass:dense_lifting}, for every
\(\psi\in W^{1,1}(\Omega)\),
\begin{equation}\label{eq:relaxation_identity_pointwise}
    J_{\rm rel}[\psi]
    =
    \overline{\mathcal J}[\psi].
\end{equation}
Consequently,
\begin{equation}\label{eq:relaxed_min_equals_unrelaxed_inf}
    \min_{\psi\in\mathcal B}\mathcal J[\psi]
    =
    \inf_{\psi\in\mathcal A}J[\psi].
\end{equation}
\end{theorem}

\begin{proof}
We first prove the lower bound.
Let \(\{\psi_k\}\subset\mathcal A\) converge to \(\psi\)
strongly in \(W^{1,1}(\Omega)\), and set
\[
    L:=\liminf_{k\to\infty}J[\psi_k].
\]
If \(L=+\infty\), then
\(\overline{\mathcal J}[\psi]\leq L\) is immediate.
Otherwise, pass to a subsequence, not relabeled, such that
\(J[\psi_k]\to L\). The total variations
\(|TD^2\psi_k|(\Omega)=\pi J[\psi_k]\) are then uniformly
bounded.

Lemma~\ref{lem:relaxed_dense_closedness} gives
\(\psi\in\mathcal B\) and
\[
    TD^2\psi_k\cws TD^2\psi
    \quad\text{in }\mathcal M(\Omega;\Rr^2),
\]
without normalizing the sequence.
Lower semicontinuity therefore yields
\[
    \overline{\mathcal J}[\psi]
    =\mathcal J[\psi]
    \leq\liminf_{k\to\infty}J[\psi_k]
    =L.
\]
Taking the infimum over all approximating sequences gives
\[
    \overline{\mathcal J}[\psi]
    \leq J_{\rm rel}[\psi].
\]

If \(\psi\notin\mathcal B\), this inequality already gives
\(J_{\rm rel}[\psi]=+\infty\).
If \(\psi\in\mathcal B\),
Theorem~\ref{thm:TD_strict_recovery} provides a sequence
\(\{\psi_k\}\subset\mathcal A\) such that
\[
    \psi_k\to\psi
    \quad\text{in }W^{1,1}(\Omega),
    \qquad
    J[\psi_k]\to\mathcal J[\psi].
\]
Hence
\[
    J_{\rm rel}[\psi]
    \leq\mathcal J[\psi]
    =\overline{\mathcal J}[\psi].
\]
This proves \eqref{eq:relaxation_identity_pointwise}.

Since \(\mathcal A\subset\mathcal B\) and
\(J=\mathcal J\) on \(\mathcal A\),
\[
    \min_{\psi\in\mathcal B}\mathcal J[\psi]
    \leq\inf_{\psi\in\mathcal A}J[\psi].
\]
Let \(\psi_*\in\mathcal B\) be a minimizer provided by
Theorem~\ref{thm:relaxed_dense_existence}.
Applying Theorem~\ref{thm:TD_strict_recovery} to
\(\psi_*\) gives a sequence
\(\{\psi_k\}\subset\mathcal A\) with
\(J[\psi_k]\to\mathcal J[\psi_*]\).
Therefore
\[
    \inf_{\psi\in\mathcal A}J[\psi]
    \leq\mathcal J[\psi_*]
    =\min_{\psi\in\mathcal B}\mathcal J[\psi].
\]
Combining the two inequalities proves
\eqref{eq:relaxed_min_equals_unrelaxed_inf}.
\end{proof}

%%%%%%%%%%%%%%%%%%%%%%%%%%%%%%%%%%%%%%%%%%%%%%%%%
\subsection{Duality and optimality conditions}
\label{sec:dense_duality}
%%%%%%%%%%%%%%%%%%%%%%%%%%%%%%%%%%%%%%%%%%%%%%%%%

We conclude the variational analysis with optimality
conditions that will be used in Section~\ref{sec:Examples}
to prove the optimality of explicit flows. The Euclidean
norm in the energy is convex but not differentiable at
zero, so we formulate these conditions through convex
duality. We first derive a dual representation of the
minimum and characterize minimizers in \(\mathcal A\).
We then obtain a sufficient condition for minimization
in the relaxed class \(\mathcal B\) when the dual field
is smooth.

Since the boundary trace is prescribed, the admissible
variations form the space
\[
    \mathcal V_0
    :=\{h\in W^{1,1}_0(\Omega):
                  TD^2h\in L^1(\Omega;\Rr^2)\},
\]
so that \(\mathcal A=\psi_b+\mathcal V_0\).
We define \(\mathcal Z\) as the set of fields
\(Z\in L^\infty(\Omega;\Rr^2)\) such that
\[
    |Z|\le1\quad\text{a.e. in }\Omega,
    \qquad
    \int_\Omega Z\cdot TD^2h\,dx\,dy=0
    \quad\text{for every }h\in\mathcal V_0.
\]
These fields provide lower bounds for the collision energy
and lead to a dual characterization of the minimum and
the minimizers in \(\mathcal A\).

\begin{proposition}[Duality and optimality]
\label{prop:dense_duality}
Under Assumptions~\ref{ass:dense_boundary}
and~\ref{ass:dense_lifting},
\[
    \min_{\psi\in\mathcal B}\mathcal J[\psi]
    =\inf_{\psi\in\mathcal A}J[\psi] =\frac1\pi\max_{Z\in\mathcal Z}
      \int_\Omega Z\cdot TD^2\psi_b\,dx\,dy.
\]
The dual maximum is attained. A function
\(\psi\in\mathcal A\) minimizes \(J\) if and only if
there is a field \(Z\in L^\infty(\Omega;\Rr^2)\) such that
\begin{equation}\label{eq:weak_EL_certificate}
\begin{gathered}
    |Z|\le1,\qquad
    Z\cdot TD^2\psi=|TD^2\psi|
    \quad\text{a.e. in }\Omega,\\
    \int_\Omega Z\cdot TD^2h\,dx\,dy=0
    \quad\text{for every }h\in\mathcal V_0.
\end{gathered}
\end{equation}

Every dual maximizer satisfies the pointwise equality in \eqref{eq:weak_EL_certificate} at every minimizer in
\(\mathcal A\).
\end{proposition}

\begin{proof}
Set \(g=TD^2\psi_b\), and let
\[
    M:=\overline{\{TD^2h:h\in\mathcal V_0\}}^{\,L^1(\Omega;\Rr^2)}.
\]
Since \(M\) is a closed linear subspace,
\[
    \pi\inf_{\psi\in\mathcal A}J[\psi]
    =\inf_{h\in\mathcal V_0}
      \|g+TD^2h\|_{L^1(\Omega)}
    =\operatorname{dist}_{L^1}(g,M).
\]
If \(g\notin M\), the Hahn--Banach theorem applied to
\(L^1(\Omega;\Rr^2)/M\) gives a bounded linear functional
on \(L^1(\Omega;\Rr^2)\) with norm one, vanishing on
\(M\), and taking the value
\(\operatorname{dist}_{L^1}(g,M)\) at \(g\).
Representing this functional by a field in
\(L^\infty(\Omega;\Rr^2)\) gives \(Z\in\mathcal Z\)
with
\[
    \int_\Omega Z\cdot g\,dx\,dy
    =\operatorname{dist}_{L^1}(g,M).
\]
If \(g\in M\), the same equality holds with \(Z=0\).
Conversely, for every \(Z\in\mathcal Z\) and
\(h\in\mathcal V_0\),
\[
    \int_\Omega Z\cdot g\,dx\,dy
    =\int_\Omega Z\cdot(g+TD^2h)\,dx\,dy
    \le\|g+TD^2h\|_{L^1(\Omega)}.
\]
Taking the infimum over \(h\in\mathcal V_0\) gives
\[
    \int_\Omega Z\cdot g\,dx\,dy
    \leq \operatorname{dist}_{L^1}(g,M)
    \quad\text{for every }Z\in\mathcal Z.
\]
The field constructed above achieves equality, so the
dual maximum is attained and
\[
    \max_{Z\in\mathcal Z}
    \int_\Omega Z\cdot g\,dx\,dy
    =\operatorname{dist}_{L^1}(g,M)
    =\pi\inf_{\psi\in\mathcal A}J[\psi].
\]
The equality with the relaxed minimum
follows from
Theorem~\ref{thm:relaxed_dense_relaxation_identity}.

For every \(\psi\in\mathcal A\) and \(Z\in\mathcal Z\),
the difference \(\psi-\psi_b\) belongs to \(\mathcal V_0\).
Hence
\begin{equation}\label{eq:dense_duality_gap}
\begin{aligned}
    &J[\psi]
      -\frac1\pi\int_\Omega Z\cdot TD^2\psi_b\,dx\,dy\\
    &\qquad=
      \frac1\pi\int_\Omega
      \bigl(|TD^2\psi|-Z\cdot TD^2\psi\bigr)\,dx\,dy
      \ge0.
\end{aligned}
\end{equation}
If \(\psi\) is a minimizer and \(Z\) is a dual maximizer,
the left-hand side vanishes. The nonnegative integrand
therefore vanishes almost everywhere, proving necessity
and the assertion concerning every dual maximizer.
Conversely, if \eqref{eq:weak_EL_certificate} holds, then
for every \(\varphi\in\mathcal A\), convexity gives
\[
    J[\varphi]-J[\psi]
    \ge\frac1\pi\int_\Omega
        Z\cdot TD^2(\varphi-\psi)\,dx\,dy
    =0.
\]
Thus \(\psi\) is a minimizer. Equation~\eqref{eq:dense_duality_gap}
also shows that the corresponding \(Z\) is a dual maximizer.
\end{proof}

The dual objective is independent of the chosen lifting:
two admissible liftings with the same trace differ by an
element of \(\mathcal V_0\).
The pointwise conditions in \eqref{eq:weak_EL_certificate}
give
\[
    Z=\frac{TD^2\psi}{|TD^2\psi|}
    \quad\text{a.e. on }\{TD^2\psi\ne0\}.
\]
On \(\{TD^2\psi=0\}\), they require only \(|Z|\le1\).

For smooth fields, the weak condition follows from an
interior equation and a natural boundary condition.
The resulting criterion also applies to measure-valued
derivatives in the relaxed class.

\begin{proposition}[Smooth sufficient condition]
\label{prop:dense_smooth_optimality}
Under Assumptions~\ref{ass:dense_boundary}
and~\ref{ass:dense_lifting}, let \(Z\) be a
\(C^2\) field on a neighborhood of \(\overline\Omega\),
with values in \(\Rr^2\), such that \(|Z|\le1\) on
\(\overline\Omega\) and
\begin{equation}\label{eq:classical_EL_boundary}
\begin{aligned}
    2\partial_{xy}Z_1+(\partial_{xx}-\partial_{yy})Z_2
    &=0 &&\text{in }\Omega,\\
    2Z_1\nu_1\nu_2+Z_2(\nu_1^2-\nu_2^2)
    &=0 &&\text{a.e. on }\Gamma.
\end{aligned}
\end{equation}
Then
\[
    \int_\Omega Z\cdot d(TD^2h)=0
\]
for every \(h\in W^{1,1}_0(\Omega)\) with
\(TD^2h\in\mathcal M(\Omega;\Rr^2)\).
If \(\psi_*\in\mathcal B\) satisfies
\[
    \int_\Omega Z\cdot d(TD^2\psi_*)
    =|TD^2\psi_*|(\Omega),
\]
then \(\psi_*\) minimizes \(\mathcal J\) on \(\mathcal B\).
If \(\psi_*\in\mathcal A\), it also minimizes \(J\)
on \(\mathcal A\).
\end{proposition}

\begin{proof}
Let \(h\in W^{1,1}_0(\Omega)\) have
\(TD^2h\in\mathcal M(\Omega;\Rr^2)\), and set
\(w=D^\perp h\). By \eqref{eq:dense_BD_norm},
\(w\in BD(\Omega)\), and \(\div w=0\).
On \(\Gamma\), \(w\) denotes its \(BD\) trace;
see \cite[Theorem~2.3]{barrosoRelaxationTheoremSpace2000}
for the trace theorem and Gauss--Green formula.
The Gauss--Green formula gives, for every
\(\eta\in C^\infty(\overline\Omega)\),
\[
\begin{aligned}
    \int_\Gamma\eta\,w\cdot\nu\,ds
    =\int_\Omega w\cdot D\eta\,dx\,dy =\int_\Omega(h_y\eta_x-h_x\eta_y)\,dx\,dy
    =0,
\end{aligned}
\]
where the last equality follows by first-order integration
by parts and the zero trace of \(h\). Therefore
\begin{equation}\label{eq:zero_trace_variation_flux}
    w\cdot\nu=0
    \quad\text{a.e. on }\Gamma.
\end{equation}

Set
\[
    S:=
    \begin{pmatrix}
        Z_1&-Z_2\\
        -Z_2&-Z_1
    \end{pmatrix}.
\]
By the definition of \(S\),
\[
    \int_\Omega Z\cdot d(TD^2h)
    =
    \int_\Omega\operatorname{tr}(S\,dEw).
\]
Applying the Gauss--Green formula and then integrating
by parts in \(h\) gives
\[
\begin{aligned}
    \int_\Omega Z\cdot d(TD^2h)
    &=\int_\Gamma(S\nu)\cdot w\,ds
      -\int_\Omega\div S\cdot D^\perp h\,dx\,dy\\
    &=\int_\Gamma(S\nu)\cdot w\,ds +
    \int_\Omega h\,
    \bigl(2\partial_{xy}Z_1
          +(\partial_{xx}-\partial_{yy})Z_2\bigr)\,dx\,dy.
\end{aligned}
\]
The interior equation in \eqref{eq:classical_EL_boundary}
makes the last integral vanish. Since
\[
    (S\nu)\cdot\tau
    =-2Z_1\nu_1\nu_2-Z_2(\nu_1^2-\nu_2^2)=0
    \quad\text{a.e. on }\Gamma,
\]
the vector \(S\nu\) is normal to \(\Gamma\).
Equation~\eqref{eq:zero_trace_variation_flux} therefore
makes the boundary integral vanish as well.

For \(\psi\in\mathcal B\), apply this identity to
\(h=\psi-\psi_*\). The bound \(|Z|\le1\) gives
\[
\begin{aligned}
    \pi\mathcal J[\psi]
    \ge\int_\Omega Z\cdot d(TD^2\psi) =\int_\Omega Z\cdot d(TD^2\psi_*)
    =\pi\mathcal J[\psi_*].
\end{aligned}
\]
Thus \(\psi_*\) minimizes the relaxed energy. If
\(\psi_*\in\mathcal A\), the equality
\(J=\mathcal J\) on \(\mathcal A\) gives the final assertion.
\end{proof}

The natural boundary condition in
\eqref{eq:classical_EL_boundary} reflects the freedom
to prescribe the normal derivative of a smooth zero-trace
variation on each smooth boundary segment.

%%%%%%%%%%%%%%%%%%%%%%%%%%%%%%%%%%%%%%%%%%%%%%%%%%%%%%%
%%%%%%%%%%%%%%%%%%% EXAMPLES %%%%%%%%%%%%%%%%%%%%%%%%%%

\section{Examples}
\label{sec:Examples}

The following examples use the optimality conditions in
Section~\ref{sec:dense_duality}. We first study a radial flow, illustrating how radial
compression generates collisions at constant density and
proving that the radial profile is an exact minimizer. We then
characterize the zero-energy minimizers for boundary traces
generated by rotations and translations, and prove uniqueness
on domains other than disks. The square example shows that
uniqueness may fail when the minimum is positive: distinct
minimizers can have different spatial distributions of local
collision intensity. 
Finally, the square--disk comparison shows how domain geometry
and the natural boundary condition determine whether a given
flow is optimal.

%%%%%%%%%%%%%%%%%%%%%%%%%%%%%%%%%%%%%%%%%%%%%%%%%

\subsection{Radial flow in an annular sector}
\label{sec:curved_examples}
%%%%%%%%%%%%%%%%%%%%%%%%%%%%%%%%%%%%%%%%%%%%%%%%%

We consider radial flow through an annular sector with
prescribed inflow on the inner arc and outflow on the outer
arc. As the flow spreads outward, incompressibility forces
the radial speed to decrease with the radius. This deceleration
creates radial compression alongside tangential expansion,
producing a positive collision intensity even at constant
density. We derive the radial velocity from incompressibility
and the boundary flux, compute its collision energy, and prove
that this radial flow minimizes the energy over the full
admissible classes, including nonradial competitors.

Let \(0<\rho<R\), \(0<\alpha<2\pi\), let $e(\theta)$
be as in \eqref{etheta}, and consider
\[
    \Omega
    =
    \{r e(\theta):\rho<r<R,\ 0<\theta<\alpha\}.
\]
This is a bounded, simply connected, nonconvex Lipschitz
domain whose boundary consists of two circular arcs and
two radial sides. For \(\alpha=\pi/2\), it is a
quarter-annular passage.

\begin{figure}[ht]
\centering
\includegraphics[width=0.75\textwidth]{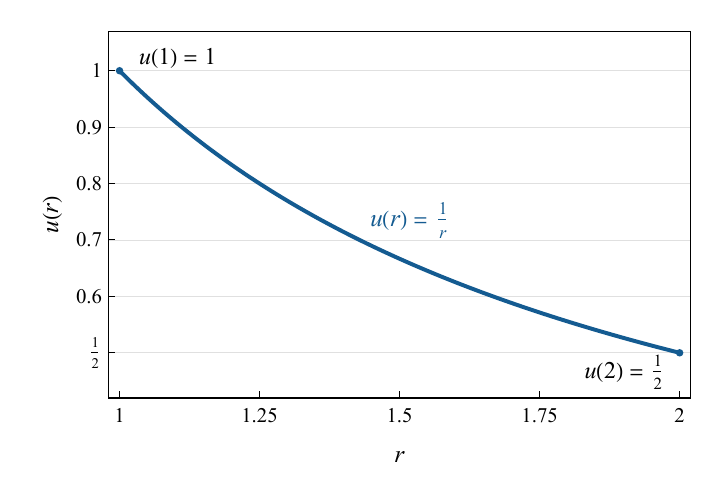}
\caption{ \small 
Radial velocity \(u(r)=1/r\) for
\(1\le r\le2\), with \(\kappa=1\).
The outward speed decreases from \(u(1)=1\) to
\(u(2)=1/2\).}
\label{fig:curved_passage_velocity}
\end{figure}

Consider a radial velocity
\[
    v(r e(\theta))=u(r)e(\theta),
    \qquad u\in C^1([\rho,R]).
\]
Its divergence is
\[
    \div v
    =
    u'(r)+\frac{u(r)}r
    =
    \frac1r(ru(r))'.
\]
Incompressibility therefore requires \(ru(r)\) to be
constant. For an outward flow, we write
\begin{equation}\label{eq:curved_radial_velocity}
    u(r)=\frac{\kappa}{r},
    \qquad
    u'(r)=-\frac{\kappa}{r^2}<0,
    \qquad \kappa>0.
\end{equation}
Thus the speed decreases from \(\kappa/\rho\) at the
inner arc to \(\kappa/R\) at the outer arc.
Figure~\ref{fig:curved_passage_velocity} shows
this decrease for \(\rho=1\), \(R=2\), and \(\kappa=1\).

{
For the radial velocity field \(v=u(r)e(\theta)\),
the matrix of \(Ev\) in the orthonormal polar frame
\((e(\theta),e'(\theta))\) is
\[
    Ev=
    \begin{pmatrix}
        u'(r)&0\\
        0&u(r)/r
    \end{pmatrix}.
\]
Its eigenvalues \(u'(r)\) and \(u(r)/r\) describe
local compression or expansion in the radial and
tangential directions, respectively.
}

By Proposition~\ref{prop:collision_integral}, the collision
functional therefore has the radial representation
\begin{equation}\label{eq:curved_radial_functional}
    \mathcal C_\Omega
    =
    \frac{\alpha}{\pi}\int_\rho^R
    r\,I\!\left(
        \diag\!\left(u'(r),\frac{u(r)}r\right)
    \right)\,dr.
\end{equation}

We prescribe the boundary flux generated by this field:
\[
\begin{aligned}
    j&=-\frac\kappa\rho
    &&\text{on the inner arc},\\
    j&=\frac\kappa R
    &&\text{on the outer arc},\\
    j&=0
    &&\text{on the radial sides}.
\end{aligned}
\]
The total inflow and outflow both equal \(\kappa\alpha\).

A boundary primitive is \(\bar j=\kappa\theta\), and
we take the lifting
\[
    \psi_b(r e(\theta))=\kappa\theta.
\]

Since \(\alpha<2\pi\) and \(\rho>0\), the function
\(\theta\) admits a smooth branch on a neighborhood of
\(\overline\Omega\). Thus
Assumption~\ref{ass:dense_lifting} holds.

Among radial velocity fields, incompressibility and the
prescribed boundary flux determine the flow uniquely. The following proposition shows
that this flow also minimizes the collision energy over
the full admissible classes.

{
\begin{proposition}[Radial flow through a curved passage]
\label{prop:radial_sector_minimum}
Let \(0<\rho<R\), \(0<\alpha<2\pi\), and \(\kappa>0\).
On
\[
    \Omega=\{r e(\theta):\rho<r<R,\ 0<\theta<\alpha\},
\]
let \(\mathcal A\) and \(\mathcal B\) be the admissible
classes \eqref{eq:unrelaxed_class_A} and
\eqref{eq:relaxed_class_B} with lifting
\(\psi_b(r e(\theta))=\kappa\theta\).
Then \(\psi_*(r e(\theta))=\kappa\theta\) attains both
minima, and
\[
    \min_{\psi\in\mathcal B}\mathcal J[\psi]
    =
    \min_{\psi\in\mathcal A}J[\psi]
    =
    \frac{2\kappa\alpha}{\pi}\log\frac R\rho.
\]
\end{proposition}
}

\begin{proof}
{
The function \(\psi_*\) is smooth on a neighborhood of
\(\overline\Omega\) and has the prescribed trace, so
\(\psi_*\in\mathcal A\). Its velocity is
\[
    D^\perp\psi_*(r e(\theta))=\frac\kappa r e(\theta).
\]
For this radial flow, the eigenvalues of \(Ev\) are
\[
    u'(r)=-\frac{\kappa}{r^2},
    \qquad
    \frac{u(r)}r=\frac{\kappa}{r^2}.
\]
Radial compression and tangential expansion balance,
preserving local area. By
Remark~\ref{rem:traceless_collision_integral}, the local
collision intensity is
\[
    \mathcal C_*(r e(\theta))
    =\frac1\pi I(Ev)
    =\frac2\pi|u'(r)|
    =\frac{2\kappa}{\pi r^2}.
\]
Integration in polar coordinates gives
\begin{equation}\label{eq:curved_radial_energy}
    J[\psi_*]
    =
    \frac{2\alpha}{\pi}\int_\rho^R r|u'(r)|\,dr
    =
    \frac{2\kappa\alpha}{\pi}\log\frac R\rho.
\end{equation}

It remains to prove minimality over the full admissible
classes.
}
Define
\[
    Z(r e(\theta))
    =(-\cos(2\theta),\sin(2\theta)).
\]
Then \(|Z|=1\), and direct differentiation gives
\[
    TD^2\psi_*=\frac{2\kappa}{r^2}Z.
\]

Since \(\rho>0\), \(Z\) is smooth on a neighborhood of
\(\overline\Omega\). Direct differentiation gives
\[
    2\partial_{xy}Z_1
    +(\partial_{xx}-\partial_{yy})Z_2=0.
\]
The unit normal is \(\nu=\pm e(\theta)\) on the circular
arcs and \(\nu=\pm(-\sin\theta,\cos\theta)\) on the radial
sides. In both cases,
\[
    2Z_1\nu_1\nu_2+Z_2(\nu_1^2-\nu_2^2)=0.
\]
Thus \(Z\) satisfies \eqref{eq:classical_EL_boundary}.
Since \(Z\cdot TD^2\psi_*=|TD^2\psi_*|\),
Proposition~\ref{prop:dense_smooth_optimality} shows that
\(\psi_*\) minimizes both \(J\) over \(\mathcal A\) and
\(\mathcal J\) over \(\mathcal B\).
\end{proof}

For \(\kappa<0\), the same conclusion holds with minimum
\[
    \frac{2|\kappa|\alpha}{\pi}\log\frac R\rho.
\]
The flow is inward, with inflow through the outer arc
and outflow through the inner arc. Indeed, negating the
stream function reverses the prescribed boundary flux
and leaves both \(J\) and \(\mathcal J\) unchanged.

\begin{figure}[ht]
\centering
\includegraphics[width=0.75\textwidth]{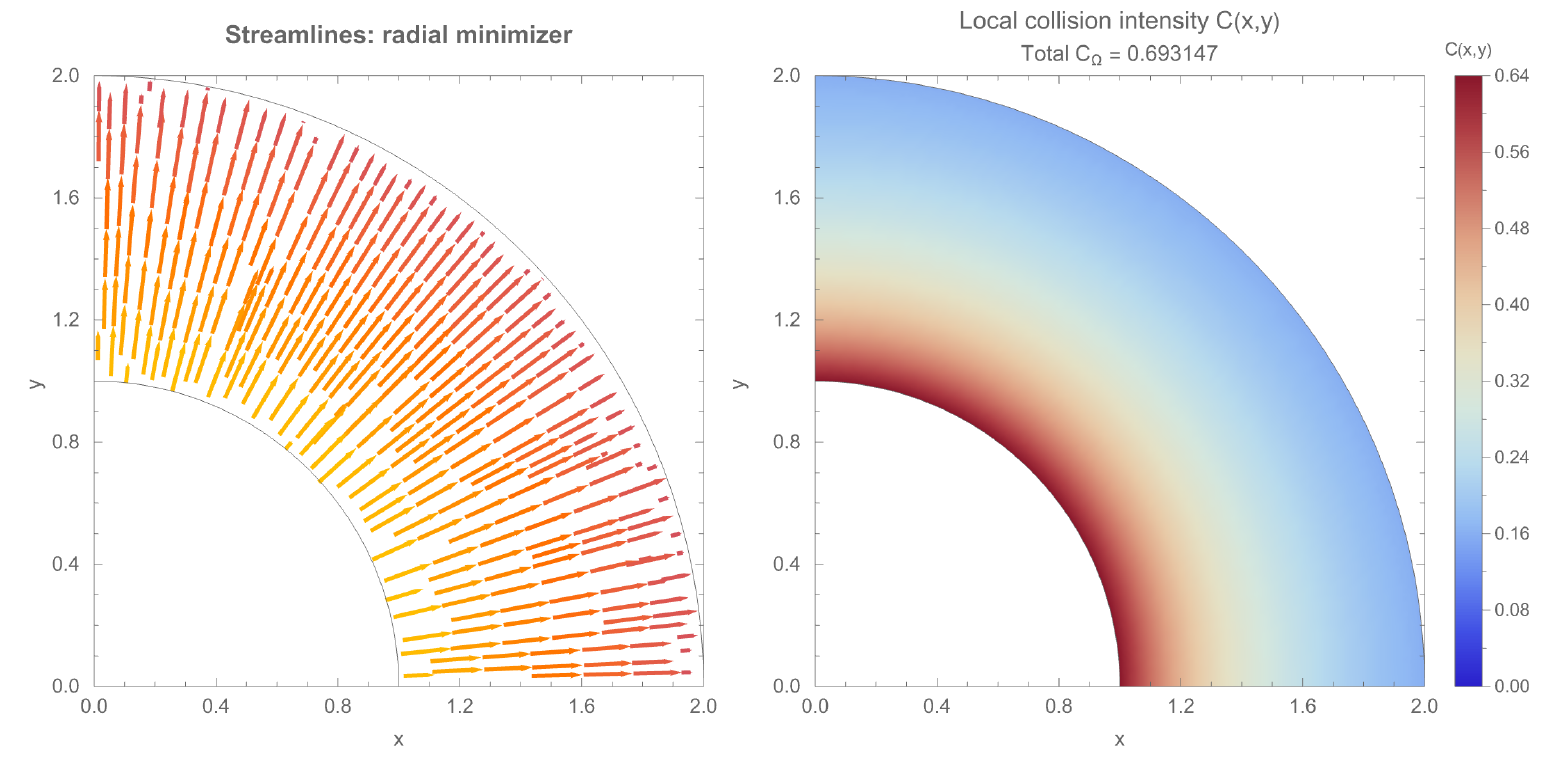}
\caption{ \small 
Radial flow through the quarter-annular passage
\(1<r<2\), \(0<\theta<\pi/2\), with \(\kappa=1\):
streamlines of \(D^\perp\psi_*\), where \(\psi_*=\theta\)
(left), and local collision intensity (right).
The flow points outward and its total collision energy is
\(J[\psi_*]=\log2\).}
\label{fig:curved_passage_radial}
\end{figure}

For \(\alpha=\pi/2\), the minimum is
\(|\kappa|\log(R/\rho)\).
Figure~\ref{fig:curved_passage_radial} shows this minimizer
for \(\rho=1\), \(R=2\), and \(\kappa=1\).
Since this domain is not a disk,
Corollary~\ref{cor:relaxed_sublevels} applies without
normalization.

\subsection{Zero collision minimizers}
\label{sec:zero_collision}

{
Rigid rotations and translations preserve distances between
particles and produce no local compression. Their symmetric
velocity gradients vanish, so their collision energy is zero
on any admissible domain. Since the energy is nonnegative,
each such flow is a minimizer for the boundary data it generates.

By Lemma~\ref{lem:dense_distributional_kernel}, every
zero-energy stream function has the form
\[
    q(x,y)=\frac{\alpha}{2}(x^2+y^2)
           +\beta x+\gamma y+\delta,
    \qquad \alpha,\beta,\gamma,\delta\in\Rr.
\]
Its velocity
\[
    D^\perp q=(\alpha y+\gamma,-\alpha x-\beta)
\]
is a rigid rotation combined with a constant translation.
The remaining question is whether the prescribed trace
allows more than one such flow. 
For the prescribed trace \(q|_\Gamma\), every minimizer
has zero energy, so uniqueness is determined by the
zero-trace kernel \(\mathcal K_0\).

\begin{corollary}[Zero-energy minimizers]
\label{cor:zero_collision_minimizers}
Let \(\Omega\subset\Rr^2\) be a bounded simply connected
Lipschitz domain, and let \(q\) be a polynomial of the form
\eqref{eq:dense_distributional_kernel}.
For the admissible classes \(\mathcal A\) and
\(\mathcal B\) with lifting \(\psi_b=q\),
\[
    \min_{\psi\in\mathcal A}J[\psi]
    =\min_{\psi\in\mathcal B}\mathcal J[\psi]=0.
\]
The minimizers in both classes are exactly the functions
in \(q+\mathcal K_0\). Thus \(q\) is the unique minimizer
when \(\Omega\) is not a disk. On a disk, all minimizers
have the form \(q+tq_0\), \(t\in\Rr\), where \(q_0\)
is given by \eqref{eq:dense_disk_kernel_generator}.
\end{corollary}

\begin{proof}
Since \(TD^2q=0\) and \(q\in\mathcal A\), both minima
are zero. A function \(\psi\in\mathcal B\) is a minimizer
if and only if \(TD^2\psi=0\). Since \(\psi-q\) has
zero trace, this is equivalent to
\(\psi-q\in\mathcal K_0\). Each such function also
belongs to \(\mathcal A\).
The conclusions for disks and other domains follow from
Proposition~\ref{prop:dense_kernel_geometry}.
\end{proof}

On a disk, the additional rotation generated by \(q_0\)
is tangent to the boundary and therefore preserves the
prescribed flux. If \(\Omega\) is not a disk, then
\(\mathcal J[\psi]>0\) for every
\(\psi\in\mathcal B\setminus\{q\}\).
}

\subsection{Minimizers on the square}
\label{sec:nonuniqueness}
%%%%%%%%%%%%%%%%%%%%%%%%%%%%%%%%%%%%%%%%%%%%%%%%%

{
The preceding corollary gives uniqueness on domains other
than disks when the prescribed trace admits zero energy.
Uniqueness need not hold when the minimum is positive.
On the square, the saddle flow belongs to a family of
distinct minimizers with the same prescribed boundary flux
and total energy, but different local collision intensities.
}

For \(\kappa>0\), the saddle stream function
\(\psi=\kappa xy\) on \(\Omega=(0,1)^2\) gives
\[
    D^\perp\psi=(\kappa x,-\kappa y),\qquad
    TD^2\psi=(2\kappa,0),\qquad
    \mathcal C\equiv\frac{2\kappa}{\pi}.
\]
Its total collision energy is \(J[\kappa xy]=2\kappa/\pi\).

For \(\kappa=1\), the following proposition applies
the optimality conditions of
Section~\ref{sec:dense_duality} to a family of distinct
minimizers with trace \(xy\).

\begin{proposition}[Non-uniqueness on the square]
\label{prop:nonuniqueness}
Let \(\Omega=(0,1)^2\), and let \(\mathcal A_{xy}\) and
\(\mathcal B_{xy}\) be the admissible classes
\eqref{eq:unrelaxed_class_A} and \eqref{eq:relaxed_class_B}
with lifting \(\psi_b=xy\). Then
\[
    \min_{\psi\in\mathcal A_{xy}}J[\psi]
    =
    \min_{\psi\in\mathcal B_{xy}}\mathcal J[\psi]
    =
    \frac2\pi.
\]
Both minima are attained by every member of the family
\begin{equation}\label{eq:psi_a_nonunique}
    \psi_a(x,y)
    =xy+a\sin(\pi x)\sin(\pi y),
    \qquad 0\le a\le\pi^{-2}.
\end{equation}
The stream functions and their local collision intensities
are pairwise distinct as \(a\) varies.
\end{proposition}

\begin{proof}
Each \(\psi_a\) is smooth and has trace \(xy\), so
\(\psi_a\in\mathcal A_{xy}\). Direct differentiation gives
\[
    TD^2\psi_a
    =
    \bigl(2+2a\pi^2\cos(\pi x)\cos(\pi y),\,0\bigr).
\]
The first component is strictly positive in the open square
for \(0\le a\le\pi^{-2}\). Hence every member has the
same field \(Z=(1,0)\), which satisfies the interior
equation in \eqref{eq:classical_EL_boundary}.
On each side of the square, the natural boundary condition
reduces to \(Z_2=0\), which also holds.

Since \(Z\cdot TD^2\psi_a=|TD^2\psi_a|\),
Proposition~\ref{prop:dense_smooth_optimality} shows that
every \(\psi_a\) minimizes both \(J\) over \(\mathcal A_{xy}\)
and \(\mathcal J\) over \(\mathcal B_{xy}\).
Direct integration gives
\[
    \mathcal J[\psi_a]=J[\psi_a]
    =
    \frac2\pi\int_\Omega
    \bigl(1+a\pi^2\cos(\pi x)\cos(\pi y)\bigr)\,dx\,dy
    =
    \frac2\pi.
\]

The local collision intensities satisfy
\[
\begin{aligned}
    \mathcal C_a(x,y)
    &=
    \frac2\pi
    \bigl(1+a\pi^2\cos(\pi x)\cos(\pi y)\bigr),\\
    \|\mathcal C_a\|_{L^\infty(\Omega)}
    &=
    \frac2\pi(1+a\pi^2).
\end{aligned}
\]
Thus distinct values of \(a\) give distinct stream functions
and local collision intensities.
\end{proof}

\begin{figure}[ht]
\centering
\includegraphics[width=0.75\textwidth]{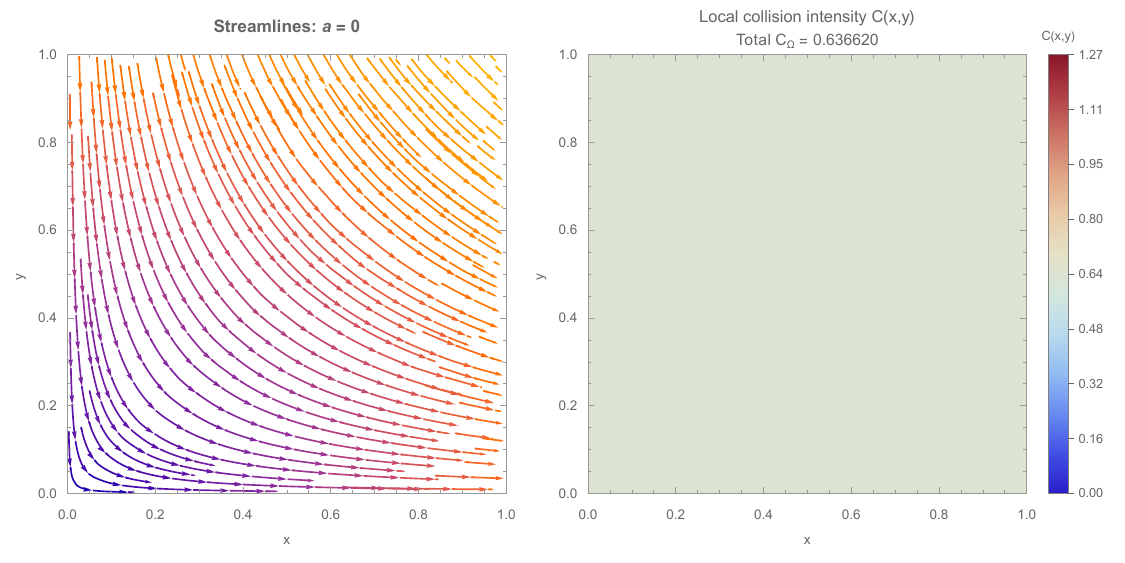}
\par\medskip
\includegraphics[width=0.75\textwidth]{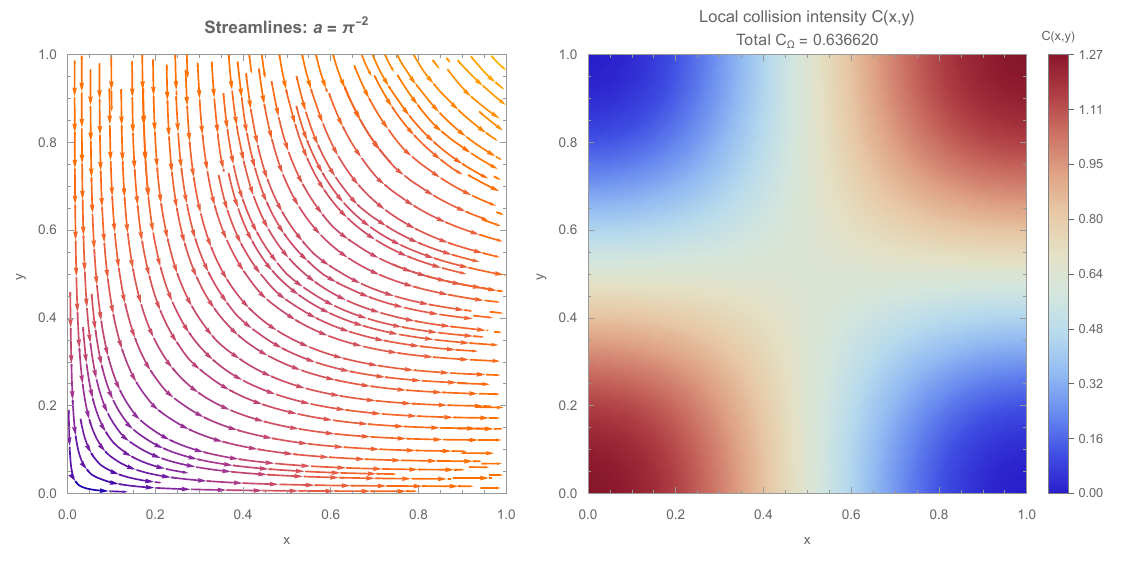}
\caption{ \small 
Two minimizers from \eqref{eq:psi_a_nonunique} with trace
\(xy\) on \((0,1)^2\): \(a=0\) (upper row) and
\(a=\pi^{-2}\) (lower row).
Streamline colors (left) use the common speed scale
\([0,\sqrt2]\), and collision intensities (right) use
\([0,4/\pi]\).
Both have energy \(2/\pi\simeq0.636620\), with intensity
suprema \(2/\pi\) and \(4/\pi\), respectively.}
\label{fig:dense}
\label{fig:nonuniqueness}
\label{fig:square_minimizer_family}
\end{figure}

By positive homogeneity, multiplying these stream functions
by \(\kappa>0\) gives minimizers for the prescribed trace
\(\kappa xy\).
By Proposition~\ref{prop:dense_kernel_geometry},
\(\mathcal K_0=\{0\}\) on the square. Thus a trivial
zero-trace kernel does not imply uniqueness of the minimizer.

Figure~\ref{fig:dense} compares the endpoints
\(a=0\) and \(a=\pi^{-2}\) of the minimizing family.
This comparison shows that a minimizer of total collision
intensity need not minimize its spatial supremum.

%%%%%%%%%%%%%%%%%%%%%%%%%%%%%%%%%%%%%%%%%%%%%%%%%

\subsection{The saddle flow on a disk}
\label{sec:disk_saddle}
%%%%%%%%%%%%%%%%%%%%%%%%%%%%%%%%%%%%%%%%%%%%%%%%%

The optimality of a flow depends on the domain and the
prescribed boundary data. The saddle flow in
Section~\ref{sec:nonuniqueness} minimizes the collision
energy on the square. On a disk, we construct a perturbation
that preserves the prescribed boundary trace and strictly
reduces the total collision energy. {
This difference in
optimality is reflected in the natural boundary condition.
}

Let
\[
    \Omega=\{(x,y)\in\Rr^2:x^2+y^2<R^2\},
    \qquad R>0,
\]
and prescribe \(\bar j=\kappa xy\) on \(\Gamma\), with
\(\kappa>0\) and lifting \(\psi_b=\kappa xy\).
The signed boundary flux is
\[
    j=\partial_\tau\bar j=\kappa R\cos(2\theta).
\]
For the saddle stream function, the pointwise conditions in
\eqref{eq:weak_EL_certificate} force \(Z=(1,0)\).
This field satisfies the interior equation in
\eqref{eq:classical_EL_boundary}.
On the disk, however,
\(\nu=e(\theta)\), so the left-hand side of the natural
boundary condition in \eqref{eq:classical_EL_boundary}
equals \(\sin(2\theta)\).
Thus the natural boundary condition fails.
The following perturbation gives an explicit decrease
of the energy.

For \(t\ge0\), define
\[
    \psi_t(x,y)
    :=\kappa xy+txy(R^2-x^2-y^2).
\]
Every \(\psi_t\) is smooth and has trace \(\bar j\), so
these functions belong to the same admissible class and
have the same prescribed boundary flux.
Direct differentiation gives
\[
    TD^2\psi_t
    =
    \bigl(2\kappa+2t(R^2-3x^2-3y^2),\,0\bigr).
\]
The first component is nonnegative when
\(0\le t\le\kappa/(2R^2)\), but changes sign for larger
\(t\). Retaining the absolute value, polar integration gives
\[
    J[\psi_t]
    =4\int_0^R\bigl|\kappa+tR^2-3tr^2\bigr|r\,dr.
\]
In particular, \(J[\psi_0]=2\kappa R^2\).

Choose \(t_c=3\kappa/(4R^2)\), giving a simple explicit
competitor with a sign change in the first component of
\(TD^2\psi_{t_c}\). Then
\[
    \mathcal C_{t_c}(r)
    =\frac{\kappa}{2\pi}
      \left|7-9\frac{r^2}{R^2}\right|,
\]
which vanishes on the circle \(r=\sqrt7\,R/3\).
With \(u=r^2/R^2\), we obtain
\[
\begin{aligned}
    J[\psi_{t_c}]
    &=\frac{\kappa R^2}{2}
      \int_0^1|7-9u|\,du\\
    &=\frac{53}{36}\kappa R^2
    =\frac{53}{72}J[\psi_0].
\end{aligned}
\]
Thus \(\psi_0=\kappa xy\) is not a minimizer for this
prescribed boundary trace on the disk.
The admissible perturbation reduces the energy by
approximately \(26.39\%\), while increasing the peak
collision intensity by \(75\%\), since
\[
    \|\mathcal C_{t_c}\|_{L^\infty(\Omega)}
    =\frac{7\kappa}{2\pi},
    \qquad
    \|\mathcal C_0\|_{L^\infty(\Omega)}
    =\frac{2\kappa}{\pi}.
\]
By Theorem~\ref{thm:relaxed_dense_relaxation_identity},
this calculation gives
\[
    \min_{\psi\in\mathcal B}\mathcal J[\psi]
    =
    \inf_{\psi\in\mathcal A}J[\psi]
    \le\frac{53}{36}\kappa R^2.
\]
The exact value is not determined by this comparison.

Note that all these competitors already lie in the normalized
subspace \(\ker P_0\) defined by \eqref{eq:dense_kernel_projection}:
since the kernel generator
\[
    q_0(x,y)=\frac{x^2+y^2-R^2}{2}
\]
is radial and \(\psi_t\) is odd in \(x\), we have
\[
    \int_\Omega\psi_tq_0\,dx\,dy=0,
    \qquad P_0\psi_t=0.
\]
Thus the energy reduction comes from a genuine interior
deformation orthogonal to the rotational kernel \(\mathcal K_0\).

\begin{figure}[ht]
\centering
\includegraphics[width=0.75\textwidth]{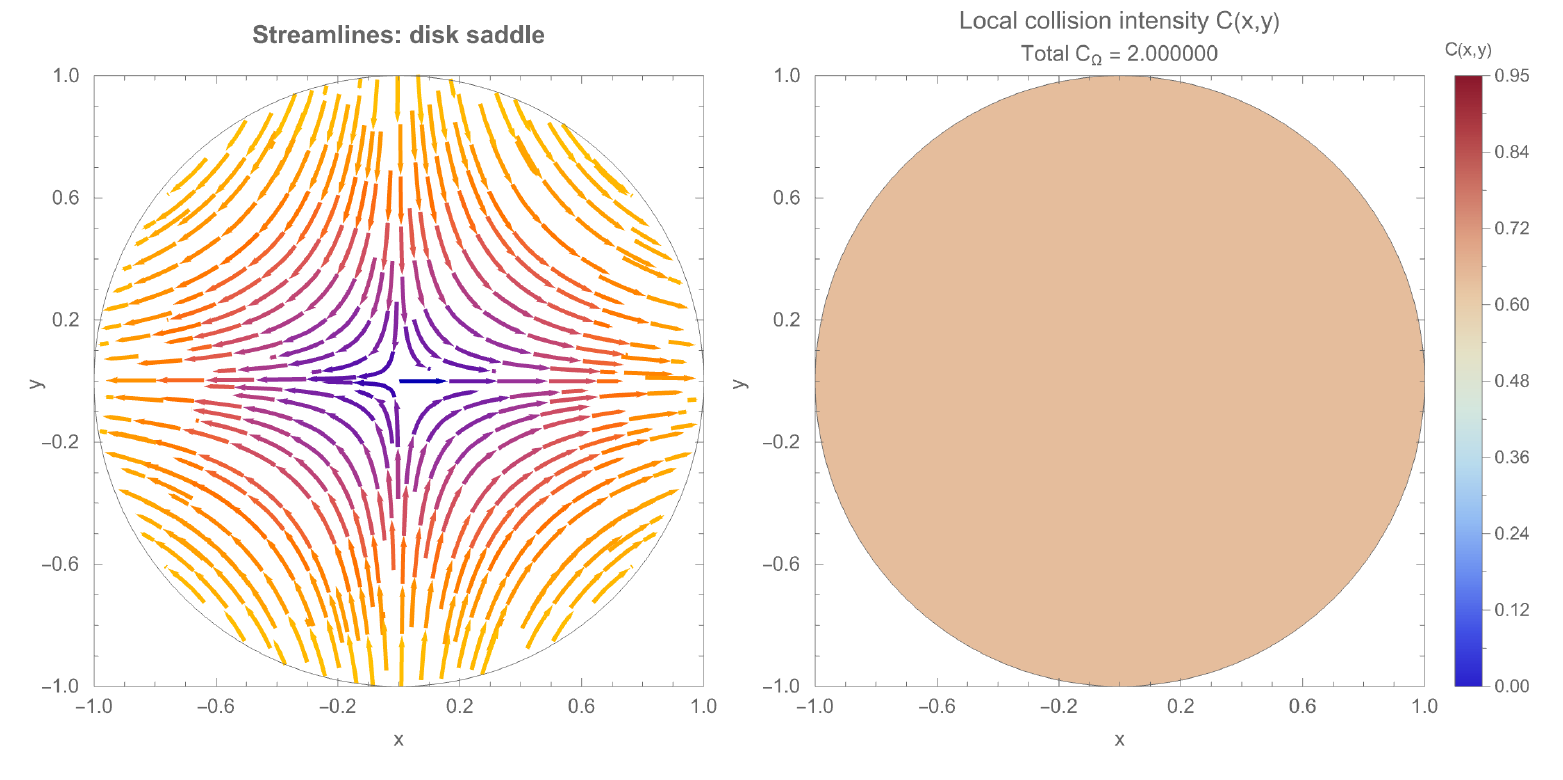}
\par\medskip
\includegraphics[width=0.75\textwidth]{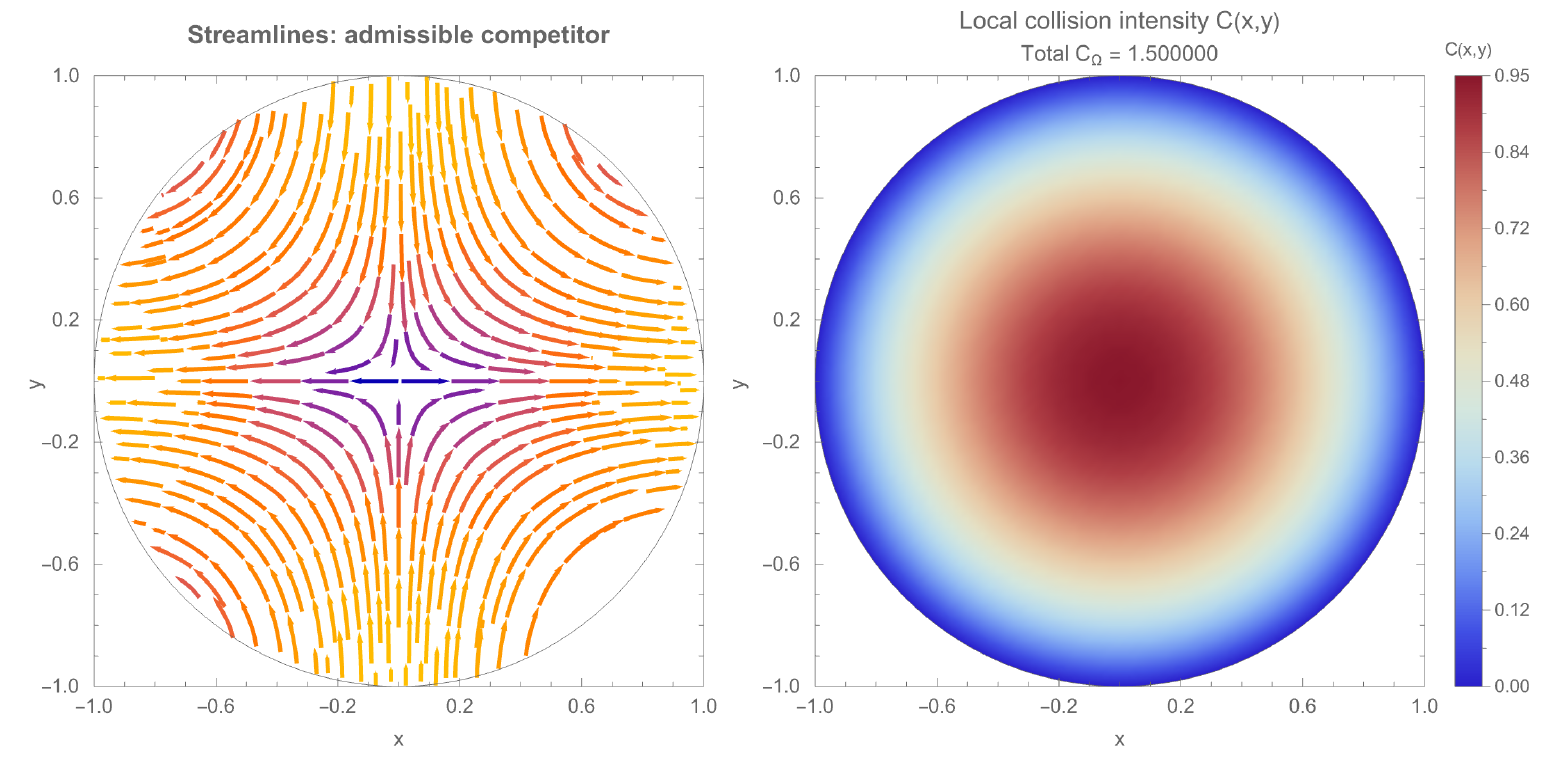}
\caption{ \small 
Saddle flow (\(t=0\), upper row) and an admissible competitor
(\(t=3/4\), lower row) on the unit disk, with \(R=\kappa=1\).
Streamlines (left) are colored by speed, with respective
ranges \([0,1]\) and \([0,7\sqrt7/18]\).
Local collision intensities (right) share the scale
\([0,7/(2\pi)]\).
At the same prescribed boundary flux, the energy decreases
from \(2\) to \(53/36\simeq1.4722\), a reduction of
approximately \(26.39\%\).}
\label{fig:disk_saddle_comparison}
\label{fig:disk_saddle}
\label{fig:disk_saddle_perturbed}
\end{figure}

Figure~\ref{fig:disk_saddle_comparison} compares \(t=0\)
and \(t=3/4=0.75\) for \(R=\kappa=1\).

%%%%%%%%%%%%%%%%%%%%%%%%%%%%%%%%%%%%%%%%%%%%%%%%%

\FloatBarrier

\section{Conclusions and future directions}
\label{sec:conclusion}

{
We introduced a collision functional for stationary flows and
studied its minimization in the dense incompressible regime
\(m\equiv1\), under prescribed boundary flux. On bounded simply
connected Lipschitz domains, we proved existence of relaxed
minimizers and exact relaxation in the strong
\(W^{1,1}(\Omega)\) topology. Coercivity and compactness follow
after removing the zero-trace kernel, with normalization
required only on disks.
The existence theorem allows concentration of the controlled
derivatives and does not establish smoothness of minimizing
velocities. For such relaxed solutions, collision intensity
is represented by a measure.

The examples show that a trivial zero-trace kernel does not
imply uniqueness. Distinct minimizers may have the same
prescribed boundary flux and total energy while producing
different spatial distributions of local collision intensity.
They also illustrate how the geometry of the domain and the
natural boundary condition affect optimality.

Natural extensions include the study of variable-density flows
in addition to the dense incompressible setting considered here.
Another direction is the optimization of domain geometry,
inflow and outflow portions, and boundary flux distributions,
for both constant and variable densities, to reduce the
minimum total collision intensity while keeping the total
inflow fixed. Shape optimization will also require geometric
constraints that prevent domain degeneration.
A mean-field-game extension would require specifying
individual costs and comparing the resulting equilibria
with the collision-minimizing flows studied here.
}

\bibliographystyle{abbrv}
%\bibliography{mfgv11_nn}
\bibliography{supbib}

\end{document}